\documentclass[aap,citesort,MSNbibl,dvips]{arximspdf}
\usepackage{mathbh}
\usepackage{accents}
\usepackage{graphicx}

\doi{10.1214/10-AAP732}
\volume{21}
\issue{6}
\pubyear{2011}
\firstpage{2075}
\lastpage{2108}

\makeatletter

\newtheorem{theorem}{Theorem}
\newtheorem{lem}{Lemma}
\newproclaim{rem}{Remark}

\newcommand{\overset}[1]{\accentset{\circ}{#1}}
\newcommand{\overlineG}{\overline{\G}{}}
\newcommand{\overlineGtwo}{\overlineG^{\hspace*{1pt}2}}

\newcommand{\llangle}{\langle\!\langle}
\newcommand{\rrangle}{\rangle\!\rangle}

\newcommand{\cyl}{\operatorname{cyl}}
\newcommand{\card}{\operatorname{card}}
\newcommand{\hyp}{\operatorname{hyp}}
\newcommand{\flow}{\operatorname{flow}}

\newcommand{\eps}{\varepsilon}
\newcommand{\PP}{\mathbb{P}}
\newcommand{\RR}{\mathbb{R}}
\newcommand{\EE}{\mathbb{E}}
\newcommand{\NN}{\mathbb{N}}
\newcommand{\ZZ}{\mathbb{Z}}
\newcommand{\E}{\mathcal{E}}
\newcommand{\N}{\mathcal{N}}
\newcommand{\I}{\mathcal{I}}
\newcommand{\G}{\Gamma}
\newcommand{\p}{\partial}
\newcommand{\C}{\mathcal{C}}
\newcommand{\V}{\mathcal{V}}
\newcommand{\U}{\mathcal{U}}
\newcommand{\M}{\mathcal{M}}

\newcommand{\W}{\mathcal{W}}

\newcommand{\wideparen}{\overbrace}
\newcommand{\Pco}{\hspace*{3pt}\overset{\hspace*{-3pt}\wideparen{{\RR}^d\setminus P}}}
\newcommand{\Po}{\overset{P}}

\newcommand{\bro}{\overset{B}}

\makeatother

\begin{document}
\begin{frontmatter}

\title{Upper large deviations for the maximal flow through a domain of
$\bolds{\mathbb{R}^d}$ in first passage percolation}
\runtitle{Upper large deviations for the maximal flow}

\begin{aug}
\author[A]{\fnms{Rapha\"{e}l} \snm{Cerf}\ead[label=e1]{rcerf@math.u-psud.fr}} and
\author[B]{\fnms{Marie} \snm{Th\'{e}ret}\corref{}\ead[label=e2]{marie.theret@univ-paris-diderot.fr}}
\runauthor{R. Cerf and M. Th\'{e}ret}
\affiliation{Universit\'{e} Paris Sud and \'{E}cole Normale Sup\'{e}rieure}
\address[A]{D\'{e}partement de Math\'{e}matiques\\
Universit\'{e} Paris Sud\\
Math\'{e}matiques b\^{a}timent 425\\
91405 Orsay Cedex\\
France\\
\printead{e1}}
\address[B]{D\'{e}partement de Math\'{e}matiques\\
\quad et Applications\\
\'{E}cole Normale Sup\'{e}rieure\\
45 rue d'Ulm\\
75230 Paris Cedex 05\\
France\\
\printead{e2}}
\end{aug}

\received{\smonth{7} \syear{2009}}
\revised{\smonth{3} \syear{2010}}

%
\begin{abstract}
We consider the standard first passage percolation model in the
rescaled graph
$\ZZ^d/n$ for $d\geq2$ and a domain $\Omega$ of boundary $\G$ in~%
$\RR^d$. Let $\G^1$ and $\G^2$ be two disjoint open subsets of $\G$ representing
the parts of $\G$ through which some water can enter and escape
from~$\Omega$. We investigate the asymptotic behavior of the flow $\phi_n$
through a~discrete version $\Omega_n$ of $\Omega$ between the
corresponding discrete
sets $\G^1_n$ and $\G^2_n$. We prove that under some
conditions on the regularity of the domain and on the law of the
capacity of
the edges, the upper large deviations of $\phi_n/ n^{d-1}$ above a
certain constant are of volume order, that is, decays exponentially
fast with
$n^d$. This article is part of a larger project\vspace*{1pt} in which the authors prove
that this constant is the a.s. limit of $\phi_n/n^{d-1}$.
\end{abstract}

\setattribute{keyword}{AMS}{AMS 2000 subject classification.}

\end{frontmatter}

\section{First definitions and main result}

We use notation introduced in \cite{KestenStFlour}
and~\cite{Kestenflows}. Let $d\geq2$. We consider the graph
$(\mathbb{Z}^{d}_n, \mathbb E ^{d}_n)$ having for vertices $\mathbb Z
^{d}_n = \ZZ^d/n$ and for edges $\mathbb E ^{d}_n$, the set of pairs of
nearest neighbors for the standard~$L^{1}$ norm. With each edge $e$ in
$\mathbb{E}^{d}_n$ we associate a random variable~$t(e)$ with values in
$\mathbb{R}^{+}$. We suppose that the family $(t(e), e
\in\mathbb{E}^{d}_n)$ is independent and identically distributed with a
common law $\Lambda$; this is the standard model of first passage
percolation on the graph $(\mathbb{Z}^d_n, \mathbb{E}^d_n)$. We
interpret $t(e)$ as the capacity of the edge $e$; it means that $t(e)$
is the maximal amount of fluid that can go through the edge $e$ per
unit of time.

\begin{figure}

\includegraphics{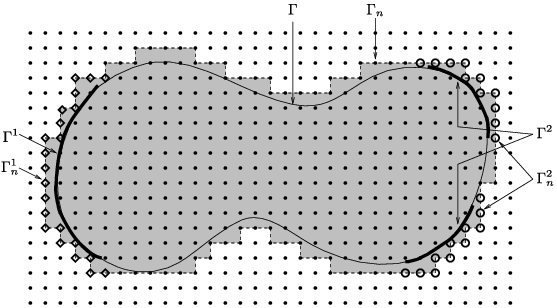}

\caption{Domain $\Omega$.}
\label{chapitre7domaine2}
\end{figure}

We consider an open bounded connected subset $\Omega$ of $\RR^d$ such that
the boundary $\G= \p\Omega$ of $\Omega$ is piecewise
of class $\C^1$ [in particular $\G$ has finite area, $\mathcal
{H}^{d-1}(\G)
<\infty$]. It means that $\G$ is included in the union of a finite
number of
hypersurfaces of class $\C^1$, that is, in the union of a finite number
of~%
$C^1$ submanifolds of $\RR^d$ of codimension $1$. Let $\G^1$, $\G^2$ be
two disjoint subsets of~$\G$ that are
open in $\G$.
We want to define the maximal flow from $\G^1$ to $\G^2$ through $\Omega
$ for
the capacities $(t(e), e\in\EE^d_n)$. We consider a discrete version
$(\Omega_n, \G_n, \G^1_n, \G^2_n)$ of $(\Omega, \G, \G^1,\G^2)$ defined by
\[
\cases{
\Omega_n = \{ x\in\ZZ^d_n \mid
d_{\infty}(x,\Omega) <1/n \} ,\cr
\G_n = \{ x\in
\Omega_n \mid  \exists y \notin\Omega_n , \langle x,y \rangle\in\EE^d_n
\} ,\cr
\G^i_n = \{ x\in\G_n \mid  d_\infty(x, \G^i) <1/n ,
d_\infty(x, \G^{3-i}) \geq1/n \}, &\quad for $i=1,2 $,}\
\]
where $d_{\infty}$ is the $L^{\infty}$-distance and the notation
$\langle
x,y\rangle$ corresponds to the edge of endpoints $x$ and $y$ (see
Figure \ref{chapitre7domaine2}).

%

We shall study the maximal flow from $\G^1_n$ to $\G^2_n$ in $\Omega_n$.
Let us properly define
the maximal flow $\phi(F_1 \rightarrow F_2$ in $C)$
[also denoted by $\phi(F_1 \rightarrow F_2$ in $C\cap\mathbb
{Z}^d/n)$] from
$F_1$ to $F_2$ in $C$ for $C \subset\mathbb{R}^d$ (or by commodity the
corresponding graph $C\cap\mathbb{Z}^d/n$). We will say that an edge
$e=\langle x,y\rangle$ belongs to a subset~$A$ of $\mathbb{R}^{d}$, which
we denote by $e\in A$ if the\vspace*{2pt} interior of the segment joining~$x$ to $y$
is included in $A$. We define
$\widetilde{\mathbb{E}}_n^{d}$ as the set of all the oriented edges,
that is,
an element $\widetilde{e}$ in $\widetilde{\mathbb{E}}_n^{d}$ is an ordered
pair of vertices which are nearest neighbors. We denote an element
$\widetilde{e} \in\widetilde{\mathbb{E}}_n^{d}$ by $\llangle
x,y \rrangle$, where $x$, $y \in\mathbb{Z}_n^{d}$ are the
endpoints of $\widetilde{e}$ and the edge is oriented from $x$ toward~$y$.
We consider the set $\mathcal{S}$ of all pairs of functions
$(g,o)$, with $g\dvtx\mathbb{E}_n^{d} \rightarrow\mathbb{R}^{+}$ and
$o\dvtx\mathbb{E}_n^{d} \rightarrow\widetilde{\mathbb{E}}_n^{d}$ such that
$o(\langle x,y\rangle) \in\{ \llangle x,y\rrangle,
\llangle y,x \rrangle\}$, satisfying the following:
\begin{longlist}[(1)]
\item[(1)] for each edge $e$ in $C$ we have
\[
0 \leq g(e) \leq t(e) ,
\]
\item[(2)] for each vertex $v$ in $C \setminus(F_1\cup F_2)$ we have
\[
\sum_{e\in C \dvtx o(e)=\llangle v,\cdot\rrangle}
g(e) = \sum_{e\in C \dvtx o(e)=\llangle\cdot,v \rrangle} g(e) ,
\]
\end{longlist}
where the notation $o(e) = \llangle v,\cdot \rrangle$
[resp., $o(e) = \llangle \cdot ,v \rrangle$] means that
there exists $y \in\mathbb{Z}_n^d$ such that $e = \langle v,y \rangle$
and $o(e) = \llangle v,y \rrangle$ [resp., $o(e) =
\llangle y,v \rrangle$].
A~couple $(g,o) \in\mathcal{S}$ is a possible stream in $C$ from
$F_1$ to $F_2$; $g(e)$ is the amount of fluid that goes through the
edge $e$ and $o(e)$ gives the direction in which the fluid goes through~$e$.
The two conditions on $(g,o)$ express only the fact that the
amount of fluid that can go through an edge is bounded by its capacity
and that there is no loss of fluid in the graph. With each possible
stream we associate the corresponding flow
\begin{eqnarray*}
\flow(g,o) & = & \sum_{ u \in F_2 , v \notin C \dvtx \langle
u,v\rangle\in\mathbb{E}_n^{d}} g(\langle u,v\rangle)
\mathbh{1}_{o(\langle u,v\rangle) = \llangle u,v \rrangle} -
g(\langle u,v\rangle) \mathbh{1}_{o(\langle u,v\rangle) = \llangle
v,u \rrangle} \\
& = &\sum_{ u \in F_1 , v \notin C \dvtx \langle
u,v\rangle\in\mathbb{E}_n^{d}} g(\langle u,v\rangle)
\mathbh{1}_{o(\langle u,v\rangle) = \llangle v,u \rrangle} -
g(\langle u,v\rangle) \mathbh{1}_{o(\langle u,v\rangle) = \llangle
u,v \rrangle} .
\end{eqnarray*}

This is the amount of fluid that crosses $C$ from $F_1$
to $F_2$ if the fluid respects the stream $(g,o)$. The two definitions
are equivalent since the stream satisfies the node law
at each vertex of $C\setminus(F_1 \cup F_2)$. The maximal flow through
$C$ from $F_1$ to $F_2$ is the supremum of this quantity over all
possible choices of streams,
\[
\phi(F_1 \rightarrow F_2 \mbox{ in }C) = \sup\{ \flow(g,o) \mid
(g,o) \in\mathcal{S} \} .
\]

We recall that we consider an open bounded
connected subset $\Omega$ of $\RR^d$ whose boundary $\G$ is piecewise of
class $\C^1$ and two disjoint open subsets $\G^1$ and $\G^2$ of~$\G$. We denote by
\[
\phi_n = \phi(\G^1_n \rightarrow\G^2_n \mbox{ in } \Omega_n)
\]
the maximal flow from $\G^1_n$ to $\G^2_n$ in $\Omega_n$. We will investigate
the asymptotic behavior of $\phi_n/n^{d-1}$ when $n$ goes to infinity.
More precisely, we will show that the upper large deviations of
$\phi_n$ above a certain constant $\widetilde{\phi}_{\Omega}$ are of
volume order. Here we state the precise theorem.
\begin{theorem}
\label{chapitre7devsup}
We suppose that $d(\G^1, \G^2)>0$, where $d$ is the Euclidean distance
between these two subsets of $\RR^d$. If the law $\Lambda$ of the
capacity of an edge admits an exponential moment,
\[
\exists\theta>0\qquad \int_{\RR^+} e^{\theta x} \,d\Lambda(x) < +\infty ,
\]
%
then there exists a finite constant
$\widetilde{\phi}_{\Omega}$ such that for all $\lambda>\widetilde{\phi}
_{\Omega}$,
\[
\limsup_{n\rightarrow\infty} \frac{1}{n^d} \log\PP[ \phi_n \geq
\lambda n^{d-1} ] < 0 .
\]
\end{theorem}

The description of $\widetilde{\phi}_{\Omega}$ will be given in Section
\ref{chapitre7deflimite}. As we will explain in Section~\ref{secajout3},
this constant is relevant in the sense that we prove in the companion
papers \cite{CerfTheret09geob} and \cite{CerfTheret09infb} that under added
geometric assumptions, $\widetilde{\phi}_{\Omega}$ is the almost sure
limit of
$\phi_n/n^{d-1}$. Theorem \ref{chapitre7devsup} is needed to prove this
a.s. convergence.
\begin{rem}
In Theorem \ref{chapitre7devsup} we need to impose that $d (\G^1, \G^2)
>0$ because
otherwise we cannot be sure that $\widetilde{\phi}_{\Omega} <\infty$,
as we
will see in Section \ref{chapitre7secfini}. Moreover, if
$d(\G^1, \G^2) =0$, the upper large deviations of $\phi_n/n^{d-1}$ may not
be of volume order (see Theorem \ref{chapitre3thmtau} in Section
\ref{secajout2} below).
\end{rem}
\begin{rem}
The large deviations we obtain are of the relevant order. Indeed, if all
the edges in $\Omega_n$ have a capacity which is abnormally big, then the
maximal flow $\phi_n$ will be abnormally big too. The probability for these
edges to have an abnormally large capacity is of order $\exp-C n^d$
for a~constant $C$ because the number of edges in $\Omega_n$ is $C' n^{d}$
for a
constant $C'$.
\end{rem}


The rest of the article is structured as follows. The background is
presented in Section \ref{secajout}; we first give some added definitions
in Section \ref{secajout1}, then we present the existing results concerning
maximal flows in first passage percolation in Section~\ref{secajout2} and
finally, we explain the role of this article in the comprehension of maximal
flows problems in Section \ref{secajout3}. The constant~$\widetilde{\phi}_{\Omega}$ is
computed in Section~\ref{chapitre7deflimite}. Section \ref{secsketch} gives
a detailed sketch of the proof of Theorem \ref{chapitre7devsup}. The rest
of the article is devoted to the proof itself.


\section{Background}
\label{secajout}


\subsection{Some definitions}
\label{secajout1}

\subsubsection{Geometric notation}

\begin{figure}

\includegraphics{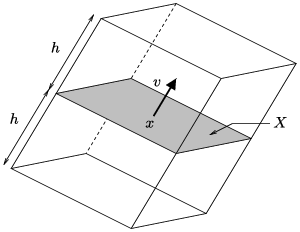}

\caption{Cylinder $\cyl(X,h)$.}
\label{chapitre7cylindrerect}
\end{figure}

We start with some geometric definitions. For a subset $X$ of
$\mathbb{R}^d$, we denote by $\mathcal{H}^s (X)$ the $s$-dimensional
Hausdorff measure of $X$ (we will use $s=d-1$ and $s=d-2$). The
$r$-neighborhood~$\V_i(X,r)$ of $X$ for the distance $d_i$ that can be
the Euclidean distance if $i=2$ or the $L^\infty$-distance if $i=\infty$,
is defined by
\[
\V_i (X,r) = \{ y\in\RR^d \mid  d_i(y,X)<r\} .
\]
If $X$ is a subset of $\RR^d$ included in a hyperplane of $\RR^d$ and of
codimension $1$ (e.g., a~nondegenerate hyperrectangle), we denote by
$\hyp(X)$ the hyperplane spanned by $X$ and we denote by $\cyl(X, h)$ the
cylinder of basis $X$ and of height $2h$ defined by
\[
\cyl(X,h) = \{x+t v \mid  x\in X , t\in
[-h,h] \} ,
\]
where $v$ is one of the two unit vectors orthogonal to $\hyp(X)$ (see
Figure \ref{chapitre7cylindrerect}).
%
For $x\in
\RR^d$, $r\geq0$ and a unit vector $v$, we denote by $B(x,r)$ the
Euclidean closed ball centered at $x$ of radius $r$.


\subsubsection{Flow in a cylinder}

Here are some particular definitions of flows through a box. It is important
to know them because all our work consists of comparing the maximal flow
$\phi_n$ in $\Omega_n$ with the maximal flows in small cylinders.
Let $A$ be a nondegenerate hyperrectangle,
that is, a box of dimension $d-1$ in $\mathbb{R}^d$. All
hyperrectangles will be
supposed to be closed in $\mathbb{R}^d$. We denote by
$v$ one of
the two unit vectors orthogonal to $\hyp(A)$. For $h$ a
positive real number, we consider the cylinder $\cyl(A,h)$.
The set $\cyl(A,h) \setminus\hyp(A)$ has two connected
components which we denote by $\mathcal{C}_1(A,h)$ and
$\mathcal{C}_2(A,h)$. For $i=1,2$, let $A_i^h$ be
the set of the points in $\mathcal{C}_i(A,h) \cap\mathbb{Z}_n^d$ which have
a~nearest neighbor in $\mathbb{Z}_n^d \setminus\cyl(A,h)$,
\[
A_i^h = \{x\in\mathcal{C}_i(A,h) \cap
\mathbb{Z}_n^d \mid  \exists y \in\mathbb{Z}_n^d \setminus\cyl(A,h)
, \langle x,y \rangle\in\EE^d_n\} .
\]
Let $T(A,h)$ [resp., $B(A,h)$] be the top
(resp., the bottom) of $\cyl(A,h)$, that is,
\begin{eqnarray*}
T(A,h) &=& \{ x\in\cyl(A,h) \mid  \exists y\notin\cyl(A,h) ,\\
&&\hspace*{3pt}
\langle x,y\rangle\in\mathbb{E}_n^d \mbox{ and }\langle x,y\rangle
\mbox{ intersects } A+hv \}
\end{eqnarray*}
and
\begin{eqnarray*}
B(A,h) &=& \{ x\in\cyl(A,h) \mid  \exists y\notin\cyl(A,h) ,\\
&&\hspace*{3pt}\langle x,y\rangle\in\mathbb{E}_n^d \mbox{ and } \langle x,y\rangle
\mbox{ intersects } A-hv \} .
\end{eqnarray*}
For a given realization $(t(e),e\in\mathbb{E}_n^{d})$, we define the variable
$\tau(A,h) = \tau(\cyl(A$, $h), v)$ by
\[
\tau(A,h) = \tau(\cyl(A,h), v) = \phi\bigl(A_1^h \rightarrow A_2^h
\mbox{ in } \cyl(A,h)\bigr)
\]
and the variable $\phi(A,h)= \phi(\cyl(A,h), v)$ by
\[
\phi(A,h) = \phi(\cyl(A,h), v) = \phi\bigl(B(A,h) \rightarrow T(A,h)
\mbox{ in } \cyl(A,h)\bigr) ,
\]
where $\phi(F_1 \rightarrow F_2$ in $C)$ is the maximal
flow from $F_1$ to $F_2$ in $C$, for $C \subset\mathbb{R}^d$ (or by
commodity the corresponding graph $C\cap\mathbb{Z}^d/n$) defined
previously. The dependence in $n$ is implicit here, in fact we can
also write $\tau_n (A,h)$ and~$\phi_n(A,h)$ if we want to emphasize
this dependence on the mesh of the graph.


\subsubsection{Max-flow min-cut theorem}

The maximal flow
$\phi(F_1\rightarrow F_2$ in $C)$ can be expressed differently
thanks to the max-flow min-cut theorem \cite{Bollobas}. We need some
definitions to state this result.
A path on the graph $\mathbb{Z}_n^{d}$ from~$v_{0}$ to $v_{m}$ is a
sequence $(v_{0}, e_{1}, v_{1},\ldots, e_{m}, v_{m})$ of vertices
$v_{0},\ldots, v_{m}$ alternating with edges $e_{1},\ldots, e_{m}$ such that
$v_{i-1}$ and $v_{i}$ are neighbors in the graph joined by the edge
$e_{i}$, for $i$ in $\{1,\ldots, m\}$.
A set $E$ of edges in $C$ is said to cut $F_1$ from $F_2$ in
$C$ if there is no path from $F_1$ to $F_2$ in $C \setminus
E$. We call $E$ an $(F_1,F_2)$-cut if $E$ cuts $F_1$ from $F_2$ in $C$
and if no proper subset of $E$ does. With each set $E$ of edges we
associate its capacity which is the variable
\[
V(E) = \sum_{e\in E} t(e) .
\]
The max-flow min-cut theorem states that
\[
\phi(F_1\rightarrow F_2\mbox{ in } C) = \min\{ V(E) \mid  E
\mbox{ is a } (F_1,F_2)\mbox{-cut} \} .
\]


\subsection{State of the art}
\label{secajout2}

\subsubsection{Existing laws of large numbers}

In this section and the next one we consider the
standard first passage percolation model on the graph $(\ZZ^d, \EE^d)$
instead of the rescaled graph $(\ZZ^d_n, \EE^d_n)$.

In dimension two, classical problems of distance in
first passage percolation and problems of flows are closely related. By the
max-flow min-cut theorem, we know that the maximal flow from the top to the
bottom of a cylinder is equal to the minimal capacity of a set of edges
that cuts the top from the bottom of the cylinder and we can notice that
in dimension two the dual of such a cutset---for the standard duality of
planar graphs---is a path from the left-hand side to the right-hand
side of the dual
cylinder; thus, if we give the same capacity at an edge and its dual but
we interpret it as the time needed to cross the dual edge in the dual
graph, the minimal capacity of a cutset in the original cylinder is equal
to the minimal time needed to go from left to right in the dual
cylinder. Thus,
results concerning maximal flows in two-dimensional first passage
percolation can be obtained thanks to the known results concerning problems
of distance in first passage percolation (see \cite{GrimmettKesten84}
and \cite{KestenStFlour}). However, such
a correspondence does not exist in
dimension three or more; the ``dual'' of an edge in dimension three is
a~small plaquette (as defined by Kesten \cite{Kestenflows}) orthogonal to
this edge that cuts it in its middle of side-length one and whose sides
are parallel to the axis of coordinates. Thus, the dual of a cutset is a
``surface'' of plaquettes and there exists no result concerning this
kind of
object in the literature of first passage percolation in terms of problems
of distance.

The existing results concerning maximal flows in first passage percolation
in dimension three or more follow. The maximal flow has been
studied almost exclusively through cylinders, since cylinders have good
properties of symmetry and stacking. The most natural flow to study in a
cylinder~$\cyl(nA,\allowbreak h(n))$ (where the height function $h$ satisfies
$\lim_{n\rightarrow\infty h(n)}=+\infty$) is $\phi(nA,h(n))$ but
$\tau(nA,h(n))$ has better properties so it is easer to deal with; it is
almost subadditive.

Using a subadditive argument and concentration inequalities, Rossignol
and Th\'{e}ret have proved in \cite{RossignolTheret08b} that $\tau(nA,
h(n))$ satisfies a law of large numbers.
\begin{theorem}[(Rossignol and Th\'{e}ret)]
We suppose that
\[
\int_{[0,\infty[} x \,d\Lambda(x) < \infty.
\]
For every unit vector $v$, for every nondegenerate hyperrectangle $A$
orthogonal to~$v$, for
every height function $h\dvtx \NN\rightarrow\RR^+$ satisfying
$\lim_{n\rightarrow\infty} h(n) = +\infty$, we have
\[
\lim_{n\rightarrow\infty} \frac{\tau(nA, h(n))}{\mathcal{H}^{d-1}(nA)}
= \nu(v) \qquad\mbox{in } L^1.
\]
Moreover, if the origin of the graph belongs to $A$ or if
\[
\int_{[0,\infty[} x^{1+{1}/({d-1})} \,d\Lambda(x) < \infty,
\]
then
\[
\lim_{n\rightarrow\infty} \frac{\tau(nA, h(n))}{\mathcal{H}^{d-1}(nA)}
= \nu(v) \qquad\mbox{a.s.}
\]
\end{theorem}

Indeed, thanks to the max-flow min-cut theorem, we know that
$\tau(nA,h(n))$ is equal to the minimal capacity of a $((nA)^{h(n)}_1,
(nA)^{h(n)}_2)$-cutset. Roughly speaking, such a cutset has its
boundary fixed along $\partial(nA)$ which lies bet\-ween $(nA)^{h(n)}_1$
and $(nA)^{h(n)}_2$. This property implies that $\tau(nA,h(n))$ is \mbox{almost}
subadditive, thus, its convergence is not surprising. In the case
where~$h(n)$ is negligible compared to $n$, then $\phi(nA,h(n))$ satisfies the
same law of large numbers as $\tau(nA,h(n))$ since the cylinder
$\cyl(nA,h(n))$ is asymptotically very flat, thus, a
$(B(nA,h(n)),T(nA,h(n)))$-cutset has also its boundary very close to
$\partial(nA)$.

We recall some geometric properties of the map $\nu\dvtx v\in S^{d-1}
\mapsto\nu(v)$ under the only condition on $\Lambda$ that $\EE
(t(e))<\infty$. They
have been stated in \cite{RossignolTheret08b}, Section 4.4. There
exists a unit vector $v_0$ such that
$\nu(v_0)=0$ if and only if for all
unit vector $v$, $\nu(v)=0$ and it happens if and only if $\Lambda(0)
\geq
1-p_c(d)$, where $p_c(d)$ is the critical parameter of the bond
percolation on $\ZZ^d$. This property has been proved by Zhang \cite
{Zhang}. Moreover, $\nu$ satisfies the weak
triangle inequality, that is, if $(ABC)$ is a nondegenerate triangle in
$\mathbb{R}^d$ and~$v_A$, $v_B$ and $v_C$ are the
exterior normal unit vectors to the sides $[BC]$, $[AC]$, $[AB]$ in the
plane spanned by $A$, $B$, $C$, then
\[
\mathcal{H}^1 ([AB]) \nu(v_C) \leq \mathcal{H}^1 ([AC])
\nu(v_B) + \mathcal{H}^1 ([BC]) \nu(v_A).
\]
This implies that the homogeneous extension $\nu_0$ of $\nu$ to $\RR^d$,
defined by $ \nu_0(0) =0 $ and for all $ w$ in $\RR^d$,
\[
\nu_0(w) = |w|_2 \nu(w/|w|_2)
\]
is a convex function; in particular, since $\nu_0$ is finite, it is
continuous on $\RR^d$. We denote by $\nu_{\min}$ (resp.,
$\nu_{\max}$) the infimum (resp., supremum) of $\nu$ on~$S^{d-1}$.

Kesten, Zhang, Rossignol and Th\'{e}ret have studied the maximal flow
between the top and the
bottom of straight cylinders. Let us denote by $D(\mathbf{k}, m)$ the cylinder
\[
D(\mathbf{k}, m) = \prod_{i=1}^{d-1} [0,k_i] \times[0,m],
\]
where\vspace*{1pt} $\mathbf{k} = (k_1,\ldots,k_{d-1}) \in\RR^{d-1}$. We denote by $\phi
(\mathbf{k},m)$ the maximal flow in $ D(\mathbf{k}, m)$ from its top $
\prod_{i=1}^{d-1} [0,k_i] \times\{ m \}$ to its bottom $
\prod_{i=1}^{d-1} [0,k_i] \times\{0\}$. Kesten \cite{Kestenflows}
proved the following result.
\begin{theorem}[(Kesten)]
Let $d=3$. We suppose that $\Lambda(0)<p_0$ for some fixed $p_0\geq
1/27$ and that
\[
\exists\gamma>0\qquad \int_{[0,+\infty[} e^{\gamma x} d \Lambda(x)
< \infty.
\]
If $m=m(\mathbf{k})$ goes to infinity with $k_1 \geq k_2$ in such a way that
\[
\exists\delta>0\qquad \lim_{k_1\geq k_2 \rightarrow\infty}
k^{-1+\delta} \log m(\mathbf{k}) = 0,
\]
then
\[
\lim_{k_1\geq k_2 \rightarrow\infty} \frac{\phi(\mathbf{k},m)}{k_1 k_2} = \nu((0,0,1)) \qquad\mbox{a.s. and in }L^1.
\]
Moreover, if $\Lambda(0)>1-p_c(d)$, where $p_c(d)$ is the critical
parameter for the
standard bond percolation model on $\ZZ^d$ and if
\[
\int_{[0,+\infty[} x^6 \,d\Lambda(x) < \infty,
\]
there exists a constant $C=C(F) <\infty$ such that for all $m=m(\mathbf{k})$
that goes to infinity with $k_1\geq k_2$ and satisfies
\[
\liminf_{k_1\geq k_2 \rightarrow\infty} \frac{m(\mathbf{k})}{k_1 k_2} > C
\]
for all $k_1\geq k_2$ sufficiently large, we have
\[
\phi(\mathbf{k},m) = 0 \qquad\mbox{a.s.}
\]
\end{theorem}

Zhang \cite{Zhang07} improved this result. 
%
\begin{theorem}[(Zhang)]
Let $d\geq2$. We suppose that
\[
\exists\gamma>0\qquad \int_{[0,+\infty[} e^{\gamma x} \,d\Lambda(x) < \infty
.
\]
Then for all $m=m(\mathbf{k})$ that goes to infinity when all the $k_i$,
$i=1,\ldots,d-1$ go to infinity in such a way that
\[
\exists\delta\in\ ]0,1]\qquad \log m(\mathbf{k})
\leq\max_{i=1,\ldots,d-1} k_i^{1-\delta},
\]
we have
\[
\lim_{k_1,\ldots,k_{d-1} \rightarrow\infty} \frac{\phi(\mathbf{k},
m)}{\prod_{i=1}^{d-1} k_i} = \nu((0,\ldots,0,1)) \qquad\mbox{a.s. and
in } L^1.
\]
Moreover, this limit is positive if and only if $\Lambda(0)<1-p_c(d)$.
\end{theorem}

To show this theorem, Zhang first obtains an important control on the
number of edges in a minimal cutset. Finally, Rossignol and Th\'{e}ret
\cite{RossignolTheret08b} improved Zhang's result in the particular
case where the dimensions of the basis of the straight cylinder go to
infinity all at the same speed. They obtain the following result.
\begin{theorem}[(Rossignol and Th\'{e}ret)]
We suppose that
\[
\int_{[0,\infty[} x \,d\Lambda(x) < \infty.
\]
For every straight hyperrectangle $A = \prod_{i=1}^{d-1} [0,a_i] \times
\{0 \}$ with $a_i >0$ for all~$i$, for every height function $h\dvtx\NN
\rightarrow
\RR^+$ satisfying $\lim_{n\rightarrow\infty} h(n) = +\infty$ and\break
$\lim_{n\rightarrow\infty} \log h(n) / n^{d-1} =0$, we have
\[
\lim_{n\rightarrow\infty} \frac{\phi(nA,h(n))}{\mathcal{H}^{d-1} (nA)}
= \nu((0,\ldots,0,1)) \qquad\mbox{a.s. and in }L^1.
\]
\end{theorem}

In dimension two more results are known, as explained previously. We
present here the two results that do not come from the literature of
problems of distance in first passage percolation. Rossignol and Th\'
{e}ret \cite{RossignolTheret09} studied the maximal flow from the
top to the bottom of a tilted cylinder in dimension two and proved the
following theorem (\cite{RossignolTheret09}, Corollary 2.10).
\begin{theorem}[(Rossignol and Th\'{e}ret)]
\label{thmd2}
Let $v$ be a unit vector, let $A$ be a~nondegenerate line-segment
orthogonal to
$v$ and $h\dvtx\NN\rightarrow\RR^+$ a height function satisfying $\lim_{n
\rightarrow\infty} h(n) = +\infty$ and
$\lim_{n\rightarrow\infty} \log h(n) / n =0$. We suppose that there
exists $\alpha\in[0,\pi/2]$ such that
\[
\lim_{n\rightarrow\infty} \frac{2h(n)}{\mathcal{H}^1(nA)} = \tan\alpha
.
\]
Then, if
\[
\int_{[0,\infty[} x \,d\Lambda(x) < \infty,
\]
we have
\[
\lim_{n \rightarrow\infty} \frac{\phi(nA,h(n))}{\mathcal{H}^1(nA)} =
\inf\biggl\{ \frac{\nu(v')}{v\cdot v'} \Bigm| v' \mbox{ satisfies
} v\cdot v' \geq\cos\alpha\biggr\} \qquad\mbox{in }L^1.
\]
Moreover, if the origin of the graph is the middle of $A$, or if
\[
\int_{[0,\infty[} x^2 \,d\Lambda(x) < \infty,
\]
then we have
\[
\lim_{n \rightarrow\infty} \frac{\phi(nA,h(n))}{\mathcal{H}^1(nA)} =
\inf\biggl\{ \frac{\nu(v')}{v\cdot v'} \Bigm| v' \mbox{ satisfies
} v\cdot v' \geq\cos\alpha\biggr\} \qquad\mbox{a.s.}
\]
\end{theorem}
%

Garet \cite{Garet2} studied the maximal flow $\sigma(A)$ between a
convex bounded set~$A$ and infinity in the case $d=2$. By an extension of
the max-flow min-cut theorem to nonfinite graphs, Garet
\cite{Garet2} proves that this maximal flow is equal to the minimal
capacity of a
set of edges that cuts all paths from $A$ to infinity. Let $\p A$ be the
boundary of $A$ and $\p^* A$ the set of the points $x \in\p A$ at which
$A$ admits a unique exterior normal unit vector $v_A(x)$ in a measure
theoretic sense (see \cite{CerfStFlour}, Section~13, for a precise
definition). If $A$ is a convex set, the set $\p^* A$ is also equal to
the set of the points $x\in\p A$ at which $A$ admits a unique exterior
normal vector in the classical sense and this vector is $v_A(x)$. Garet proved
the following theorem.
\begin{theorem}[(Garet)]
\label{thmgaret}
Let $d=2$. We suppose that $\Lambda(0) <1-p_c(2) = 1/2$ and that
\[
\exists\gamma>0\qquad \int_{[0,+\infty[} e^{\gamma x} \,d\Lambda
(x) < \infty.
\]
Then for all convex bounded set $A$ containing $0$ in its interior, we have
\[
\lim_{n\rightarrow\infty} \frac{\sigma(nA)}{n} = \int_{\p^* A}
\nu(v_A(x)) \,d\mathcal{H}^1 (x) = \I(A) > 0 \qquad\mbox{a.s.}
\]
Moreover, for all $\eps>0$ there exist constants $C_1$, $C_2>0$ depending
on $\eps$ and~$\Lambda$ such that
\[
\forall n\geq0\qquad \PP\biggl[\frac{\sigma(nA)}{n\I(A)} \notin
\ ]1-\eps, 1+\eps[ \biggr] \leq C_1 \exp(-C_2 n).
\]
\end{theorem}

Nevertheless, a law of large numbers for the maximal flow from the top
to the bottom of a tilted cylinder for $d\geq3$ was not proved yet. In
fact, the lack of symmetry of
the graph induced by the slope of the box is a major issue to extend
the existing results
concerning straight cylinders to tilted cylinders. The theorem of Garet was
not extended to dimension $d\geq3$ either.


\subsubsection{Large deviations results}

The upper and lower large deviations of the maximal flows $\phi
(nA,h(n))$ and
$\tau(nA,h(n))$ have been studied in the cases where the laws of large
numbers are known. Let us start with the existing results concerning the
upper large deviations of these maximal flows. They are studied in
\cite{Theretuppertau} and \cite{TheretUpper}. Theorem 4 in \cite
{Theretuppertau} deals with the upper large deviations
of the variable $\phi(nA,h(n))$ above $\nu(v)$.
\begin{theorem}
\label{thmdevsupphi}
We suppose that
\[
\exists\gamma>0\qquad \int_{[0,+\infty[} e^{\gamma x} \,d\Lambda(x) < \infty.
\]
Then for every unit
vector $v$ and every nondegenerate hyperrectangle $A$ orthogonal
to $v$, for every height function $h\dvtx \NN\rightarrow\RR^+$ such that
$\lim_{n\rightarrow\infty} h(n) = +\infty$ and for every
$\lambda> \nu(v)$ we have
\[
\liminf_{n\rightarrow\infty} \frac{-1}{ \mathcal{H}^{d-1}(nA) h(n) }
\log\PP
\biggl[\frac{\phi(nA,h(n))}{\mathcal{H}^{d-1}(nA)} \geq\lambda
\biggr] > 0.
\]
\end{theorem}

We shall rely on this result for proving Theorem
\ref{chapitre7devsup}. Moreover, Theorem \ref{chapitre7devsup} is a
generalization of Theorem \ref{thmdevsupphi} where we work in the domain
$\Omega$ instead of a parallelepiped. We stress the fact that $\nu(v)$
is not
in general the a.s. limit on $\phi(nA,h(n))/\mathcal{H}^{d-1}(nA)$. The
corresponding large deviation principle is proved only in the case of
straight cylinders (see Theorems 2 and 3 in \cite{TheretUpper}
\textit{that are gathered here}).
\begin{theorem}
We consider the maximal flow $\phi^{[h]}(n)$ through the straight
cylinder $[0,n]^{d-1} \times[0,h(n)]$ from its top to its bottom. We
suppose that the height function $h\dvtx\NN\rightarrow\RR^+$ satisfies
\[
\lim_{n \rightarrow\infty} \frac{h(n)}{\log n} = +\infty.
\]
Then for every $\lambda$ in $\RR^+$, the limit
\[
\psi(\lambda) = \lim_{n\rightarrow\infty} \frac{-1}{n^{d-1}h(n)}
\log\PP\bigl[\phi^{[h]}(n) \geq\lambda n^{d-1}\bigr]
\]
exists and is independent of $h$. Moreover, $\psi$ is convex on $\RR^+$,
finite and continuous on the set $\{\lambda \mid  \Lambda
([\lambda,+\infty[)>0 \}$. If
\[
\int_{[0,+\infty[} x \,d\Lambda(x) < \infty,
\]
then $\psi$ vanishes on $[0,\nu((0,\ldots,0,1))]$. If
\[
\exists\gamma>0\qquad \int_{[0,+\infty[} e^{\gamma x} \,d\Lambda(x)
< \infty,
\]
then $\psi$ is strictly positive on $]\nu((0,\ldots,0,1)),+\infty[$, and the
sequence
\[
\biggl( \frac{\phi^{[h]}(n)}{n^{d-1}} \biggr)_{n\in\NN}
\]
satisfies a large deviation principle with speed $n^{d-1}h(n)$ and governed
by the good rate function $\psi$.
\end{theorem}

The lack of symmetry makes it difficult to extend this large deviation
principle to the case of tilted cylinders. The upper large deviations
for $\tau$ depend a lot on the
moments of $\Lambda$, as proved in Theorem 3 of \cite{Theretuppertau}.
\begin{theorem}
\label{chapitre3thmtau}
Let $A$ be a nondegenerate hyperrectangle and $\vec{v}$ one of the two
unit vectors normal to $A$. Let $h\dvtx\mathbb{N} \rightarrow\mathbb{R}^+$
be a height function satisfying $\lim_{n \rightarrow \infty} h(n) =
+\infty$. The upper large
deviations of $\tau(nA,h(n))/\mathcal{H}^{d-1}(nA)$ depend on the
tail of the distribution of the capacities. Indeed, we obtain that:

\begin{longlist}
\item
If the law of the capacity of the edges has bounded support, then
for every
$\lambda>\nu(\vec{v})$ we have
\[
\liminf_{n\rightarrow\infty} \frac{-1}{\mathcal{H}^{d-1}(nA) \min
(h(n), n)} \log
\mathbb{P} \biggl[ \frac{\tau(nA, h(n))}{\mathcal{H}^{d-1}(nA)} \geq\lambda
\biggr] > 0 ;
\]
the upper large deviations are then of volume order for height functions
$h$ such that $h(n)/n$ is bounded and of order $n^d$ if
$\lim_{n\rightarrow\infty} h(n)/n = +\infty$.

\item If the capacity of the edges follows the exponential law of
parameter~$1$, then there exists $n_0 (d, A,h)$ and for
every $\lambda>\nu(\vec{v})$ there exists a positive constant $D$
depending only on $d$ and $\lambda$ such that for all $n\geq n_0$ we have
\[
\mathbb{P} [ \tau(nA, h(n)) \geq\lambda\mathcal{H}^{d-1}(nA)
] \geq \exp(- D \mathcal{H}^{d-1}(nA)).
\]

\item If the law of the capacity of the
edges satisfies
\[
\forall\theta>0 \qquad\int_{[0,+\infty[} e^{\theta x} \,d\Lambda(x) < \infty,
\]
then for all $\lambda> \nu(\vec{v})$ we have
\[
\lim_{n\rightarrow\infty} \frac{1}{\mathcal{H}^{d-1}(nA)} \log
\mathbb{P} \biggl[ \frac{\tau(nA, h(n))}{\mathcal{H}^{d-1}(nA)} \geq\lambda
\biggr] = -\infty.
\]
\end{longlist}
\end{theorem}

This dependence on the moment conditions on $\Lambda$ comes from the fact
that the distance between the part of the boundary of the cylinder through
which the water can enter and the part through which it can escape is null.

We just say a few words about lower large deviations, since it is
not the purpose of this article. Rossignol and Th\'{e}ret
\cite{RossignolTheret08b} proved that under some moment conditions, the lower
large deviations of $\tau(nA,h(n))/\mathcal{H}^{d-1}(nA)$
[resp., for
$\phi(nA,h(n))/\mathcal{H}^{d-1}(nA)$ when $h(n)$ is negligible
compared to $n$ or when $\cyl(nA,h(n))$ is straight] are of surface order,
that is, the probability that these rescaled flows are abnormally small decays
exponentially fast with $\mathcal{H}^{d-1}(nA)$. They also prove the
corresponding
large deviation principles.

From now on, we work in the rescaled graph $(\ZZ^d_n, \EE^d_n)$ again.


\vspace*{-3pt}\subsection{Global project}
\label{secajout3}

This article is part of a global project in which we prove that
$\phi_n/n^{d-1}$ converges a.s. toward $\widetilde{\phi}_{\Omega}$ and we
find the right order of the upper and lower large deviations. Indeed, we
prove in \cite{CerfTheret09infb} the following result.\vspace*{-3pt}
\begin{theorem}
\label{chapitre7devinf}
If the law $\Lambda$ of the capacity of an edge admits an exponential moment,
\[
\exists\theta>0 \qquad\int_{\RR^+} e^{\theta x} \,d\Lambda(x) < +\infty
\]
and if $\Lambda(0) < 1-p_c(d)$, then there exists a finite constant
$\phi_{\Omega}$ such that for all $\lambda< \phi_{\Omega}$,
\[
\limsup_{n\rightarrow\infty} \frac{1}{n^{d-1}}\log\PP[ \phi_n \leq
\lambda n^{d-1} ] < 0.\vspace*{-3pt}
\]
\end{theorem}

The definition of the constant $\phi_{\Omega}$ is given in
\cite{CerfTheret09infb}; it is of the same kind as the one of
$\widetilde{\phi}_{\Omega}$ we will give in Section \ref{chapitre7deflimite}
but slightly different. Finally we prove in \cite{CerfTheret09geob}
through a completely geometrical study that the constants $\phi_{\Omega
}$ and
$\widetilde{\phi}_{\Omega}$ are equal and we investigate when they are
strictly positive. Thus, by a simple Borel--Cantelli's lemma, we obtain in
\cite{CerfTheret09geob} the following result.\vspace*{-3pt}
\begin{theorem}
\label{chapitre7lgn}
We suppose that $\Omega$ is a Lipschitz domain and that $\G$ is
included in the union of a
finite number of oriented hypersurfaces $\mathcal{S}_1,\ldots,\mathcal
{S}_r$ of class $\C^1$
which are transverse to each other. We also\vspace*{2pt} suppose that $\G^1$ and~$\G^2$\vadjust{\goodbreak}
are open in $\G$, that their relative boundaries $\p_{\G} \G^1$ and
$\p_{\G}\G^2$ in $\G$ have null $\mathcal{H}^{d-1}$ measure and that
$d(\G^1,
\G^2)>0$.
We suppose that the law $\Lambda$ of the capacity of an edge admits an
exponential moment
\[
\exists\theta>0\qquad \int_{\RR^+} e^{\theta x} \,d\Lambda(x) < +\infty.
\]
Then there exists a finite constant $\phi_{\Omega} \geq0$ such that
\[
\lim_{n\rightarrow\infty} \frac{\phi_n}{n^{d-1}} = \phi_{\Omega}
\qquad\mbox{a.s.}
\]
Moreover, this equivalence holds:
\[
\phi_{\Omega} > 0 \quad\iff\quad \Lambda(0)<1-p_c(d).
\]
\end{theorem}

Combining Theorems \ref{chapitre7lgn}, \ref{chapitre7devinf} and
\ref{chapitre7devsup}, we prove that the rescaled maximal
flow $\phi_n/n^{d-1}$ converges a.s. toward a constant $\phi_{\Omega
}$, that
its upper large deviations are of volume order and that its lower large
deviations are of surface order. These theorems apply to the maximal
flow from the top to the bottom of a tilted
cylinder. Thus, they generalize the existing results concerning the
variable $\phi(A,h)$ in straight cylinders, in the particular case
where all the dimensions of the cylinder go to infinity at the same
speed (or, equivalently, the cylinder is fixed and the mesh of the
graph go to zero isotropically). The large deviation principles from above
and below still remain to be proved.


\vspace*{-3pt}\section{\texorpdfstring{Computation of $\widetilde{\phi}_{\Omega}$}
{Computation of phi Omega}}
\label{chapitre7deflimite}

\begin{figure}[t]

\includegraphics{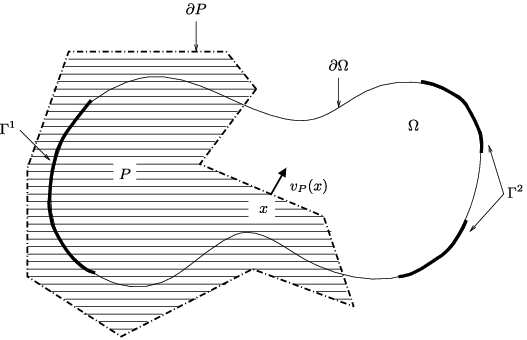}

\caption{A polyhedral set $P$ as in the definition of $\widetilde{\phi}
_{\Omega}$.}
\label{chapitre7separantsup}
\vspace*{-3pt}
\end{figure}

We give here a definition of $\widetilde{\phi}_{\Omega}$ in
terms of the map~$\nu$. When a hypersurface $\mathcal{S}$ is piecewise of
class $\C^1$, we say that $\mathcal{S}$ is
transverse to $\G$ if for all $x\in\mathcal{S} \cap\G$, the
normal unit vectors to $\mathcal{S}$ and $\G$ at $x$ are not collinear;
if the
normal vector to $\mathcal{S}$ (resp., to $\G$) at $x$ is not
well defined, this property must
be satisfied by all the vectors which are limits of normal unit vectors to
$\mathcal{S}$ (resp., $\G$) at $y \in\mathcal{S}$
(resp., $y\in\G$) when we
send $y$ to $x$---there is at most a finite number of such limits. We say
that a subset $P$
of $\RR^d$ is polyhedral if its boundary $\p P $ is included in the
union of a finite number of hyperplanes. For each point $x$ of
such a set $P$ which is on the interior of one face of $\p P$, we
denote by
$v_{P}(x)$ the exterior unit vector orthogonal to $P$ at $x$. For $A
\subset\RR^d$, we denote by $\overset{A}$ the interior of $A$. We define
$\widetilde{\phi}_{\Omega}$ by
\begin{eqnarray*}
\widetilde{\phi}_{\Omega} &=& \inf\{ \I_{\Omega} (P) \mid
P \subset\RR^d, \overlineG^1 \subset
\Po, \overlineGtwo \subset\Pco,
P\mbox{ is polyhedral},\\[-2pt]
&&\hspace*{175.4pt} \p P \mbox{ is transverse to
}\G
\},
\end{eqnarray*}
where
\[
\I_{\Omega}(P) = \int_{\partial P \cap\Omega}\nu(v_P(x)) \,d\mathcal{H}^{d-1}(x).
\]
See Figure \ref{chapitre7separantsup} for an example of such a
polyhedral set $P$.

%

The definition of the constant $\widetilde{\phi}_{\Omega}$
is not very intuitive. We propose to define the notion of a continuous\vadjust{\goodbreak}
cutset to have a better understanding of this constant. We say that
$\mathcal{S}
\subset\RR^d$ cuts $\G^1$ from $\G^2$ in $\overline{\Omega}$ if
every continuous path from $\G^1$ to $\G^2$ in $\overline{\Omega}$ intersects
$\mathcal{S}$. In fact, if $P$ is a polyhedral set of $\RR^d$ such that
\[
\overlineG^1 \subset \Po \quad\mbox{and}\quad
\overlineGtwo \subset \Pco,
\]
then $\p P \cap\overline{\Omega}$ is a continuous cutset from $\G^1$
to $\G^2$
in $\overline{\Omega}$. Since $\nu(v)$ is the average
amount of fluid that can cross a
hypersurface of area one in the direction $v$ per
unit of time, it can be interpreted as the capacity of a unitary
hypersurface orthogonal to~$v$. Thus, $\I_{\Omega}(P)$ can be
interpreted as the capacity of the
continuous cutset $\p P \cap\overline{\Omega}$ defined by $P$. The
constant $\widetilde{\phi}_{\Omega}$ is the solution of a min-cut problem
because it is equal to the infimum of the capacity of a~continuous
cutset that
satisfies some specific properties.

We remark that the capacity $\I_{\Omega}$ of a continuous cutset is
exactly the
same as the one defined by Garet in \cite{Garet2} in dimension two (see
Theorem \ref{thmgaret}), except
that we consider a maximal flow through a bounded domain so our capacity
is adapted to the problems of boundaries that arise. Moreover,
$\widetilde{\phi}_{\Omega}$ has the same form as the limit observed in
dimension two in Theorem \ref{thmd2}.


\vspace*{-3pt}\section{Sketch of the proof}
\label{secsketch}

To prove Theorem \ref{chapitre7devsup}, we have to study the
probability
%
\begin{equation}
\label{proba4}
\PP[ \phi_n \geq(\widetilde{\phi}_{\Omega} + \eps) n^{d-1}]
\end{equation}
for a positive $\eps$.\vspace*{6pt}

\textit{Step} 1. We first prove that $\widetilde{\phi}_{\Omega}$ is
finite, that is, that there exists a
polyhedral set $P \subset\RR^d$ such that $\p P$
is transverse to $\G$ and
\[
\overlineG^1 \subset \Po,\qquad \overlineGtwo \subset
\Pco.\vadjust{\goodbreak}
\]
For that purpose, at each point $x$ of $\overlineG^1$ we associate a~%
cube of center $x$, of strictly positive side-length, which is transverse
to $\G^1$ and at positive distance of $\G^2$ [this is possible thanks to
the regularity of $\G$ and the fact that $d(\G^1,\G^2)>0$]. From this
family of cubes we extract by compactness a~finite family of cubes that
covers $\overlineG^1$. The set $P$ which is defined as the union of the
cubes satisfies the desired properties.\vspace*{6pt}

\textit{Step} 2. We consider a polyhedral set $P$ as in the
definition of~$\widetilde{\phi}_{\Omega}$ such that $\I_{\Omega} (P)$
is very
close to this constant. We want to construct sets of edges near~%
$\partial P \cap\Omega$ that cut $\Gamma^1_n$ from $\Gamma^2_n$ in
$\Omega_n$. Because we took a discrete approximation of $\Omega$ from the
outside, we need to enlarge a little $\Omega$, because some flow might go
from $\G^1_n$ to $\G^2_n$ using paths that lies partly in $\Omega_n
\setminus\Omega$. Thus, we construct a set $\Omega'$ which contains a small
neighborhood of $\Omega$ (hence, also $\Omega_n$ for all $n$ large
enough) which is transverse to $\partial P$ and which is small
enough to ensure that $\I_{\Omega'} (P)$ is still very close to $\phi
_{\Omega}$. To
construct this set, we cover $\partial\Omega$ with small cubes,
by compactness we extract a finite subcover of $\partial\Omega$ and
finally we
add the cubes of the subcover to $\Omega$ to obtain $\Omega'$. We
construct these
cubes so that their boundaries are transverse to $\partial P$ and
their diameters are uniformly smaller than a small constant, so that
$\Omega'$ is included in a~neighborhood of $\Omega$ as small as we
need. Since $\partial P$ is transverse to $\G$, if we take this
constant small enough, we can control $\mathcal{H}^{d-1} (\partial P
\cap(\Omega'
\setminus\Omega))$ and thus the difference between $\I_{\Omega'} (P)$
and $\I_{\Omega} (P)$.\vspace*{6pt}

\textit{Step} 3. Then we construct a family of $C n$ (where $C>0$)
disjoint sets of edges that
cut $\Gamma^1_n$ from $\Gamma^2_n$ in $\Omega_n$ and that lie near
$\partial P$. We consider the neighborhood $P'$ of $P$ inside $\Omega'$
at distance smaller than a tiny constant~$h$ and we partition $P'
\setminus P$ into
slabs $\M'(k)$ of width of order $1/n$, so we have $ Cn$ such slabs
which look
like translations of $\partial P \cap\Omega'$ that are slightly deformed
and thickened. We prove that each
path from $\G^1_n$ to $\G^2_n$ in $\Omega_n$ must contain at least one
edge that lies in the set $\M'(k)$ for each $k$, that is, each set
$\M'(k)$ contains a cutset. Thus, we have found a family of $Cn$
disjoint cutsets.\vspace*{6pt}

\textit{Step} 4. We almost cover $\partial P \cap\Omega'$ by a finite family
of disjoint cylinders~$B_{i,j}$, whose bases are hyperrectangles
of side length $l$, that are orthogonal to $\partial P$,
of height $h$ and such that the part of $\partial P$
which is missing in this covering is very small. Thus, we obtain that
%
\begin{equation}
\label{step2beq1}
\I_{\Omega'} (P) \mbox{ is close to } \sum\nu(v_{i,j}) l^{d-1},
\end{equation}
where $v_{i,j}$ gives the direction toward which the cylinder $B_{i,j}$ is
tilted (it is the unit vector which is orthogonal to the face of
$\partial P$ that cuts $B_{i,j}$).

We want to compare $\phi_n$ with the sum of the maximal flows $\phi
(B_{i,j}, v_{i,j})$. For each $(i,j)$, let $E_{i,j}$ be a set of edges
that cuts the
top from the bottom of $B_{i,j}$. The set $\bigcup E_{i,j}$ does not cut
$\G^1_n$ from $\G^2_n$ in $\Omega_n$ in general; to create such a cutset
we must add two sets of edges.
\begin{longlist}
\item A set of edges that covers the part of $\partial P \cap
\Omega'$ that is missing in the covering by the cylinders $B_{i,j}$.
\item A set of edges that glues together all the previous sets of
edges [the sets $E_{i,j}$ and the set described in (i)].
\end{longlist}
In fact, we have already constructed $Cn$ possible sets of edges as in
(i): the edges that lie in
$\M'(k)\setminus(\bigcup B_{i,j})$ for $k=1,\ldots,Cn$. We
denote these sets by $M(k)$. We can also find $C' n$ ($C'>0$)
disjoint sets of edges that can be the glue described
in~(ii); we denote these sets by $W(l)$ for $l=1,\ldots,C'n$. Indeed, we can
choose different sets because we provide the glue more or
less in the interior of the cylinders~$B_{i,j}$. Thus, we obtain that
\begin{eqnarray}
\forall k\in\{1,\ldots,Cn\}\ \forall l\in\{1,\ldots,C'n\}\qquad
\bigcup E_{i,j} \cup M(k) \cup W(l)
\nonumber\\
&&\eqntext{\mbox{cuts }\G^1_n \mbox{ from } \G^2_n \mbox{ in
}\Omega_n.}
\end{eqnarray}
Then
%
\begin{equation}
\label{step2beq2}
\phi_n \leq \sum\phi(B_{i,j}, v_{i,j}) + \min_{k=1,\ldots,Cn}
V(M(k)) + \min_{l=1,\ldots,C'n} V(W(l)).
\end{equation}
Combining (\ref{step2beq1}) and (\ref{step2beq2}) we see that if
$\phi_n \geq(\widetilde{\phi}_{\Omega} +\eps) n^{d-1}$, one of the
following events must happen:
\begin{longlist}[(a)]
\item[(a)] $\exists j\in J$ $\phi(B_j, v_j) \geq (\nu(v_j) +
\eps/2) l^{d-1} n^{d-1} $,
\item[(b)] $\forall k\in\{1,\ldots,Cn \}$ $V(M(k)) \geq \eta
n^{d-1} $,
\item[(c)] $\forall l\in\{1,\ldots,C'n \}$ $V(W(l)) \geq \eta
n^{d-1} $,
\end{longlist}
where $\eta$ is a very small constant (depending on $\eps$ and $\phi
_{\Omega}$).\vspace*{6pt}

\textit{Step} 5. It consists in taking care of the
probability that the events (a), (b) or (c) happen. The probability of
(a) has already been studied in \cite{Theretuppertau}; cf. Theorem~\ref
{thmdevsupphi} above; the upper
large deviations of the variable $\phi$ in a cylinder above $\nu$ are
of volume order. The events (b) and (c) are of the same type
and their probability is of the form
%
\begin{equation}
\label{step3beq1}
\PP\Biggl[ \sum_{m=1}^{\alpha n^{d-1}} t_m \geq\eta n^{d-1} \Biggr]^{Dn},
\end{equation}
where $(t_m)_{m\in\NN}$ is a family of i.i.d. variables of
distribution function $\Lambda$, $D$ is a constant, $\eta$ is a very
small constant and $\alpha n^{d-1}$ is the cardinality of the family
of variables we consider. If $\alpha< \eta\EE[t_1]^{-1}$ and if the
law $\Lambda$ admits one exponential moment, the
Cram\'{e}r theorem in $\RR$ states that the probability~%
(\ref{step3beq1}) decays exponentially fast with $n^d$. Note the
role of the optimization over~$Dn$ different probabilities to
obtain the correct speed of decay. The proof would have been slightly
simpler if we would have proven only that the decay of the
probability (\ref{proba4}) is at least exponential in $n^{d-1}$.\vadjust{\goodbreak}

\textit{Step} 6.
To complete the proof, it is enough to
control the cardinality of the sets $M(k)$ and $W(l)$ for each $k$,
$l$, to ensure that we can use the Cram\'{e}r theorem as explained in
Step 5. This can be done using the geometrical properties of
$\partial P$ (it is polyhedral and transverse to $\partial\Omega'$).


\section{\texorpdfstring{The constant $\widetilde{\phi}_{\Omega}$ is finite}
{The constant phi Omega is finite}}
\label{chapitre7secfini}

To prove that $\widetilde{\phi}_{\Omega} <\infty$, it is sufficient to
exhibit a set
$P$ satisfying all the conditions given in the definition of
$\widetilde{\phi}_{\Omega}$. Indeed, if such a set $P$ exists, then
\[
\widetilde{\phi}_{\Omega} \leq \nu_{\max} \mathcal{H}^{d-1}(\p P \cap
\Omega) <
\infty
\]
since a polyhedral set has finite perimeter in $\Omega$. We will
construct such
a~set $P$. The idea of the proof is the following. We will cover
$\overlineG^1$ with small hypercubes which are transverse to $\G^1$ and
at positive distance of $\overlineGtwo$. Then, by compactness, we will
extract a finite covering. We will denote by $P$ the union of the
hypercubes of this finite covering. Then $P$
satisfies the desired properties.

We prove a geometric lemma.

\begin{lem}
Let $\Gamma$ be an hypersurface
(i.e., a $C^1$ submanifold of $\RR^d$ of codimension $1$)
and let $K$ be a compact subset of $\Gamma$.
There exists a positive $M=M(\Gamma,K)$ such that
\[
\forall\eps>0\ \exists r>0\
\forall x,y\in K\quad
|x-y|_2\leq r \quad\Rightarrow\quad
d_2(y,\tan(\Gamma,x))\leq M \eps|x-y|_2
\]
$[\tan(\Gamma,x)$ is the tangent hyperplane of $\Gamma$ at $x]$.
\end{lem}
\begin{pf}
By a standard compactness argument, it is enough to prove
the following local property:
\begin{eqnarray*}
&&
\forall x\in\Gamma\mbox{ }
\exists M(x)>0 \mbox{ }\forall\eps>0\mbox{ }
\exists r(x,\eps)>0\mbox{ }
\forall y,z\in\Gamma\cap B(x,r(x,\eps))
\\
&&\qquad d_2(y,\tan(\Gamma,z))\leq M(x) \eps|y-z|_2.
\end{eqnarray*}
Indeed, if this property holds, we cover $K$
by the open balls $\bro(x,r(x,\eps)/2)$, \mbox{$x\in K$},
we extract a finite subcover
$\bro(x_i,r(x_i,\eps)/2)$, $1\leq i\leq k$,
and we set
\[
M=\max\{ M(x_i)\dvtx1\leq i\leq k \},\qquad
r=\min\{ r(x_i,\eps)/2\dvtx1\leq i\leq k \}.
\]
Now let $y,z$ belong to $K$ with $|y-z|_2\leq r$.
Let $i$ be such that $y$ belongs to
$B(x_i,r(x_i,\eps)/2)$.
Since $r \leq r(x_i,\eps)/2$,
then both $y,z$ belong to the ball
$B(x_i,r(x_i,\eps))$ and it follows that
\[
d_2(y,\tan(\Gamma,z)) \leq M(x_i) \eps|y-z|_2\leq
M \eps|y-z|_2.
\]


We turn now to the proof of the above local property.
Since $\Gamma$ is an hypersurface, for any $x$ in $\Gamma$ there
exists a neighborhood $V$ of $x$ in $\RR^d$, a diffeomorphism
$f\dvtx V\mapsto\RR^d$ of class $C^1$ and a $(d-1)$-dimensional
vector space~$Z$ of~$\RR^d$ such that
$Z\cap f(V)=f(\Gamma\cap V)$ (see, e.g., \cite{FED}, $3.1.19$).
Let~$A$ be a~compact neighborhood of $x$ included in $V$.
Since $f$ is a diffeomorphism, the maps
$y\in A\mapsto df(y)\in\operatorname{End}(\RR^d)$,
$u\in f(A)\mapsto df^{-1}(u)\in\operatorname{End}(\RR^d)$
are continuous.
Therefore, they are bounded
\[
\exists M>0 \mbox{ }\forall y\in A \qquad\|df(y)\|\leq M,\qquad
\forall u\in f(A) \qquad\|df^{-1}(u)\|\leq M
\]
[here
$\|df(x)\|=\sup\{ |df(x)(y)|_2\dvtx|y|_2\leq1 \}$
is the standard operator norm in $\operatorname{End}(\RR^d)$].
Since $f(A)$ is compact,
the differential map $df^{-1}$ is uniformly continuous on $f(A)$,
\[
\forall\eps>0 \mbox{ }\exists\delta>0\mbox{ }
\forall u,v\in f(A)\qquad
|u-v|_2\leq\delta\quad\Rightarrow\quad
\|df^{-1}(u)-df^{-1}(v)\|\leq\eps.
\]
Let $\eps$ be positive\vspace*{1pt} and let $\delta$ be associated to $\eps$
as above.
Let $\rho$ be positive and small enough so that
$\rho<\delta/2$ and $B(f(x),\rho)\subset f(A)$
[since $f$ is a $C^1$ diffeomorphism, $f(A)$ is a
neighborhood of $f(x)$].
Let $r$ be such that $0<r<\rho/M$ and
$B(x,r)\subset A$. We claim that $M$ associated to $x$ and
$r$ associated to $\eps,x$ answer the problem.
Let $y,z$ belong to $\Gamma\cap B(x,r)$.
Since $[y,z]\subset B(x,r)\subset A$,
and $\|df(\zeta)\|\leq M$ on~$A$, then
\begin{eqnarray*}
|f(y)-f(x)|_2&\leq& M|y-x|_2\leq Mr<\rho,\qquad
|f(z)-f(x)|_2<\rho, \\
|f(y)-f(z)|_2&<&\delta,\qquad
|f(y)-f(z)|_2<M|y-z|_2.
\end{eqnarray*}
%
We next apply a classical lemma
of differential calculus (see \cite{LA}, I, 4, Corollary $2$)
to the map $f^{-1}$ and the
interval $[f(z),f(y)]$ [which is included in $B(f(x),\rho)\subset f(A)$]
and the point $f(z)$
\begin{eqnarray*}
&&\bigl|y-z-df^{-1}(f(z))\bigl(f(y)-f(z)\bigr)\bigr|_2 \\
&&\qquad\leq|f(y)-f(z)|_2
\sup\{ \|df^{-1}(\zeta)-df^{-1}(f(z))\|\dvtx\zeta\in[f(z),f(y)] \}.
\end{eqnarray*}
The right-hand member is less than
$M|y-z|_2 \eps$.
Since
$z+df^{-1}(f(z))(f(y)-f(z))$
belongs to $\tan(\Gamma,z)$, the proof is complete.
\end{pf}

We come back to our case. The boundary $\G$ of $\Omega$ is piecewise of class
$\C^1$, that is, it is included in a finite union of $\C^1$ hypersurfaces
which we denote by $(S_1,\ldots,S_p)$.
The hypersurfaces $S_1,\ldots,S_p$ being $\C^1$ and the set $\G$
compact, the maps $x\in\G\mapsto v_{S_k}(x)$, $1\leq k\leq p$
[where $v_{S_k}(x)$ is the unit normal vector to $S_k$ at $x$]
are uniformly continuous:
\begin{eqnarray*}
&&\forall\delta>0 \mbox{ }\exists\eta>0\mbox{ }
\forall
k\in\{ 1,\ldots,p \}\mbox{ }
\forall x,y\in S_k\cap\G\\
&&\qquad|x-y|_2\leq\eta\quad\Rightarrow\quad
|v_{S_k}(x)-v_{S_k}(y)|_2<\delta.
\end{eqnarray*}
Let $\eta^*$ be associated to $\delta=1$ by this property. Let
$k\in\{ 1,\ldots,p \}$. The set $S_k\cap\G$ is a compact subset of the
hypersurface $S_k$. Applying the previous lemma, we get
%
\begin{eqnarray*}
&&\exists M_k\mbox{ }
\forall\delta_0>0 \mbox{ }\exists\eta_k>0\mbox{ }
\forall x,y\in
S_k\cap\G\\
&&\qquad|x-y|_2\leq\eta_k \quad\Rightarrow\quad
d_2(y,\tan(S_k,x))
\leq M_k\delta_0 |x-y|_2.
\end{eqnarray*}
Let $M_0 = \max_{1\leq k\leq p}M_k$ and
let $\delta_0$ in $]0,1/2[$
be such that $M_0\delta_0< 1/2$.
For each $k$ in $\{ 1,\ldots,p \}$, let
$\eta_k$ be associated to
$\delta_0$ as in the above property and let
\[
\eta_0 = \min\biggl(\min_{1\leq k\leq p}\eta_k,
\eta^*,
\frac{1}{8d} \operatorname{dist}
(\G^1,\G^2)
\biggr).
\]
We build a family of cubes $Q(x,r)$, indexed by $x\in\G$
and $r\in\ ]0,r_\G[$ such that
$Q(x,r)$ is a cube centered at $x$ of side length $r$ which is transverse
to~$\G$.
For $x\in\RR^d$ and $k\in\{ 1,\ldots,p \}$, let $p_k(x)$ be a
point of $S_k\cap\G$ such that
\[
|x-p_k(x)|_2 =
\inf\{ |x-y|_2\dvtx y\in S_k\cap\G\}.
\]
Such a point exists since
$S_k\cap\G$ is compact.
We define then for
$k\in\{ 1,\ldots,p \}$
\[
\forall x\in\RR^d\qquad
v_k(x) = v_{S_k}(p_k(x)).
\]
We also define
\[
d_r =
\inf_{v_1,\ldots,v_p\in S^{d-1}}
\max_{b\in{\mathcal B}_d}
\mathop{\min_{1\leq k\leq p}}_{e\in b}
(|e-v_k|_2,|{-e-v_k}|_2),
\]
where
${\mathcal B}_d$
is the collection of the orthonormal basis of ${\mathbb R}^d$
and $S^{d-1}$ is the
unit sphere of $\RR^d$.
Let $\eta$ be associated to $d_r/4$ as in the
above continuity property. We set
\[
r_\G=\frac{\eta}{2d}.
\]
Let $x\in\G$. By the definition of $d_r$, there exists
an orthonormal basis $b_x$ of~$\RR^d$ such that
\[
\forall e\in b_x\mbox{ }
\forall k\in\{ 1,\ldots,p \}\qquad
\min\bigl(|e-v_k(x)|_2,
|{-e-v_k(x)}|_2\bigr) > \frac{d_r}{2}.
\]
Let $Q(x,r)$ be the cube centered at $x$ of sidelength $r$
whose sides are parallel to the vectors of $b_x$.
We claim that $Q(x,r)$ is transverse to $\G$ for $r<r_\G$.
Indeed, let $y\in Q(x,r)\cap\G$. Suppose that $y\in S_k$ for some
$k\in\{ 1,\ldots,p \}$, so that
$v_k(y)=v_{S_k}(y)$ and
$|x-p_k(x)|_2<dr_\G$. In particular, we have
$|y-p_k(x)|_2<2dr_\G<\eta$ and
$|v_{S_k}(y)-v_k(x)|_2<d_r/4$.
For $e\in b_x$,
\[
\frac{d_r}{2} \leq
|e-v_k(x)|_2
\leq
|e-v_{S_k}(y)|_2+
|v_{S_k}(y)-v_k(x)|_2
\]
whence,
\[
|e-v_{S_k}(y)|_2
\geq
\frac{d_r}{2}-
\frac{d_r}{4} =
\frac{d_r}{4}
.
\]
%
This is also true for $-e$, therefore,
the faces of the cube $Q(x,r)$ are transverse to~$S_k$.

Now we consider the collection
\[
\bigl(\overset{Q}(x,r), x\in\overlineG^1, r< r_{\G}\bigr).\vadjust{\goodbreak}
\]
It covers $\overlineG^1$. By compactness of
$\overlineG^1$, we can extract a finite covering $(\overset{Q}(x_i,
r_i),\allowbreak i\in I)$ from this collection. We define
\[
P = \bigcup_{i\in I} Q(x_i,r_i).
\]
We claim that $P$ satisfies all the hypotheses in the
definition of
$\widetilde{\phi}_{\Omega}$. Indeed,~$P$ is obviously polyhedral and transverse
to $\G$. Moreover, we know that
\[
\overlineG^1 \subset \Po
\]
and since $d(P, \overlineGtwo) >0$, we also obtain that
\[
\overlineGtwo \subset \Pco.
\]


\section{\texorpdfstring{Definition of the set $\Omega'$}{Definition of the set
Omega'}}

Let $\lambda$ be in $]\widetilde{\phi}_{\Omega}, +\infty[$. We are studying
\[
\PP[\phi_n \geq\lambda n^{d-1}].
\]
Suppose first that $\widetilde{\phi}_{\Omega}>0$.
There exists a positive $s$ such that $\lambda> \widetilde{\phi}_{\Omega}
(1+s)^2$. By definition of $\widetilde{\phi}_{\Omega}$, for every positive
$s$, there exists
a polyhedral subset $P$ of $\RR^d$, such that $\p P$ is transverse
to $\G$,
\[
\overlineG^1 \subset \Po,\qquad \overlineGtwo \subset \Pco
\]
and
\[
\I_{\Omega} (P) \leq \widetilde{\phi}_{\Omega} (1+s).
\]
Then $ \lambda> \I_{\Omega} (P) (1+s) $ and
\[
\PP[\phi_n \geq\lambda n^{d-1}] \leq \PP[\phi_n \geq\I_{\Omega} (P)
(1+s) n^{d-1}].
\]
Since $\p P$ is
transverse to $\G$, we know that there exists $\delta_0 >0$ (depending
on~$\lambda$, $P$ and $\G$) such that for all $\delta\leq\delta_0$,
\[
\mathcal{H}^{d-1} \bigl(\p P \cap\bigl(\V_2(\Omega, \delta) \setminus\Omega\bigr)\bigr)
\leq
\frac{s \I_{\Omega}(P)}{2 \nu_{\max}}.
\]
Thus, for any set $\Omega'$ satisfying $\Omega\subset
\Omega' \subset\V_2(\Omega, \delta_0)$, we have
\[
\int_{\p P \cap\Omega' } \nu(v_{P} (x)) \,d\mathcal{H}^{d-1}(x) \leq \I
_{\Omega}(P)
(1+ s/2),
\]
then $\lambda>( 1+s/2)( \int_{\p P \cap\Omega'} \nu(v_P (x)) \,d\mathcal
{H}^{d-1}(x))$ and
\[
\PP[\phi_n \geq\lambda n^{d-1}] \leq \PP\biggl[\phi_n \geq\biggl(
\int_{\p P \cap\Omega'} \nu(v_P (x)) \,d\mathcal{H}^{d-1}(x) \biggr) (1+s/2) n^{d-1}
\biggr].
\]

Suppose now that $\widetilde{\phi}_{\Omega}=0$. Then for an arbitrarily
fixed $s\in\ ]0,1[$, there exists a polyhedral subset $P$ of $\RR^d$,
such that $\p P$ is transverse
to $\G$,
\[
\overlineG^1 \subset \Po, \qquad\overlineGtwo \subset \Pco
\]
and
\[
\I_{\Omega} (P) \leq \frac{\lambda}{1+s}
\]
and, thus, $ \lambda> \I_{\Omega} (P) (1+s) $. If $\I_{\Omega} (P) >0$,
we can
use exactly the same argument as previously. We suppose that $\I_{\Omega}
(P)=0$. We know as previously that there exists $\delta_0 >0$
(depending on
$\lambda$, $P$ and $\G$) such that for all $\delta\leq\delta_0$,
\[
\mathcal{H}^{d-1} \bigl(\p P \cap\bigl(\V_2(\Omega, \delta) \setminus\Omega\bigr)\bigr) <
\frac{\lambda}{ \nu_{\max} (1+s/2)}.
\]
Thus, in any case, we obtain that there exists $\delta_0>0$ such that,
for any set~$\Omega'$ satisfying $\Omega\subset\Omega' \subset\V
_2(\Omega,
\delta_0)$, we have
\[
\PP[\phi_n \geq\lambda n^{d-1}] \leq \PP\biggl[\phi_n \geq\biggl(
\int_{\p P \cap\Omega'} \nu(v_P (x)) \,d\mathcal{H}^{d-1}(x) \biggr) (1+s/2) n^{d-1}
\biggr].
\]

We will construct a particular set $\Omega'$ satisfying $\Omega\subset
\Omega' \subset\V_2(\Omega, \delta_0)$. In the previous section we
have associated to each couple $(x, r)$ in $\G
\times\ ]0, r_{\G}[$ a~hypercube $Q(x,r)$ centered at $x$, of sidelength
$r$, and which is transverse to~$\G$. Using exactly the same method, we can
build a family of hypercubes\looseness=1
\[
\bigl(Q'(x,r), x\in\G, r<r_{(\G,P)} \bigr)
\]\looseness=0
such that $Q'(x,r)$ is centered at $x$, of sidelength $r$ and it is
transverse to $\G$ and $\p P$.
The family
\[
\bigl(\overset{Q}'(x,r), x\in\G, r< \min\bigl(r_{(\G,P)}, \delta_0/(2d)\bigr) \bigr)
\]
is a covering of the compact set $\G$, thus we can extract a finite covering
from this collection; we denote it by $(\overset{Q}'(x_i, r_i), i
\in J)$. We define
\[
\Omega' = \Omega\cup\bigcup_{i\in J} \overset{Q}'(x_i, r_i).
\]
Since $r_i \leq\delta_0/(2d)$ for all $i\in J$, we have $ \Omega'
\subset
\V_2(\Omega, \delta_0)$. Moreover, $\p P$ is transverse to the boundary
$\G'$ of $ \Omega'$. Finally, if
we define
\[
\delta_1 = \min_{i\in J} r_i/2,
\]
we know that $ \V_2 (\Omega,\delta_1 ) \subset\Omega' $ and thus, for
all $n\geq2d/\delta_1$, we have $\Omega_n \subset\Omega'$.

\section{\texorpdfstring{Existence of a family of $(\G^1_n,\G^2_n)$-cuts}
{Existence of a family of (Gamma^1_n,Gamma^2_n)-cuts}}

\begin{figure}

\includegraphics{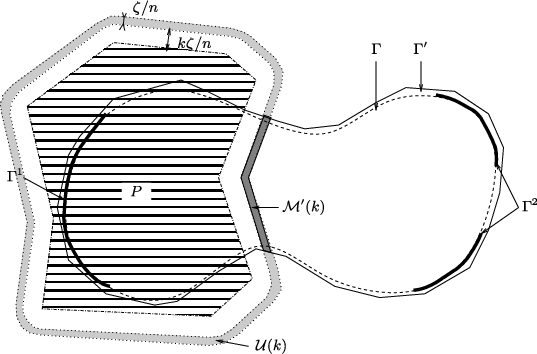}

\caption{The sets $P$, $\U(k)$ and $\M'(k)$.}
\label{chapitre7ensMbis}
\end{figure}

In this section we prove that we can construct a family of disjoint
$(\G^1_n, \G^2_n)$-cuts in $\Omega_n$. Let $\zeta$ be a fixed constant larger
than $2d$. We consider a parameter $h < h_0 = d(\p P, \G^1 \cup
\G^2)$. For $k\in\{0,\ldots,\lfloor hn/\zeta\rfloor\}$ we define
\[
P(k) = \{ x\in\RR^d \mid  d(x,P) \leq k \zeta/n \}
\]
and for $k\in\{0,\ldots,\lfloor hn/\zeta \rfloor-1 \}$ we define
\begin{eqnarray*}
\U(k) & = &(\hspace*{3pt}\overset{\hspace*{-3pt}\wideparen{{\RR
}^d\setminus P_{k+1}}})
\setminus\Po_k \\
& = &\{ x\in\RR^d \mid  k\zeta/ n \leq d(x,P) < (k+1)\zeta/n \}
\end{eqnarray*}
and $\M'(k) = \U(k)\cap\Omega'$ (see Figure \ref{chapitre7ensMbis}).
%
We will prove the following lemma.
\begin{lem}
\label{chapitre7lemcutset}
There exists $N$ large enough such that for all $n\geq N$,
every path on the graph $(\ZZ^d_n, \EE^d_n)$ from $\G^1_n$ to $\G^2_n$
in $\Omega_n$ contains at least one
edge which is included in the set $\M'(k)$ for $k\in\{0,\ldots,\lfloor
hn/\zeta
\rfloor-1 \}$.
\end{lem}

This lemma states precisely that for all $k \in
\{0,\ldots,\lfloor hn/\zeta \rfloor-1 \}$, $\M'(k)$ contains a $(\G^1_n,
\G^2_n)$-cut in $\Omega_n$.
\begin{pf*}{Proof of Lemma \ref{chapitre7lemcutset}}
Let $k\in\{0,\ldots,\lfloor hn/\zeta\rfloor-1 \}$.
Let $\gamma$ be a discrete path from $\G^1_n$ to
$\G^2_n$ in $\Omega_n$. In particular, $\gamma$ is continuous, so we can
parametrize it: $\gamma= (\gamma_t)_{0\leq t \leq1}$. There exists $N$
large enough such that for all $n\geq N$, we have
\[
\Omega_n \subset \Omega',\qquad \G^1_n \subset \V_2(\G^1,2d/n) \subset \Po
_k
\]
and
\[
\G^2_n \subset \V_2(\G^2, 2d/n) \subset
\hspace*{3pt}\overset{\hspace*{-3pt}\wideparen{{\RR}^d\setminus
P_{k+1}}}.
\]
Since $\gamma$ is continuous, we know that there exists $t_1, t_2\in\ ]0,
1[$ such that
\begin{eqnarray*}
t_1 &=& \sup\{ t\in[0,1] \mid  \gamma_t \in\Po_k \},
\\
t_2 &=& \inf\{ t\geq t_1 \mid  \gamma_t \in
\hspace*{3pt}\overset{\hspace*{-3pt}\wideparen{{\RR}^d\setminus
P_{k+1}}}\}.
\end{eqnarray*}
Since
\[
\Po_k \cup\U(k) \cup\hspace*{3pt}\overset{\hspace*{-3pt}\wideparen
{{\RR}^d\setminus
P_{k+1}}}
\]
is a partition of $\RR^d$, we know that $(\gamma_t)_{t_1 \leq t < t_2}$,
which is a continuous path, is
included in $\U(k)$. The length of $(\gamma_t)_{t_1\leq t < t_2}$ is larger
than $d(\gamma_{t_1}, \gamma_{t_2})$. The segment
$[\gamma_{t_1},\gamma_{t_2}]$ intersects
\[
\{ x\in\RR^d \mid  d(x, P) = (k+1/2)\zeta/n \}
\]
at a point $z$ and we know that
\[
\V_2\bigl(z, \zeta/(2n)\bigr) \subset \hspace*{3pt}\overset{\hspace
*{-3pt}\wideparen{V(k)}}.
\]
Thus, $d(\gamma_{t_1}, \gamma_{t_2}) \geq\zeta/n$ and then the length of
$(\gamma_t)_{t_1 \leq t < t_2}$ is larger than $\zeta/n$. Finally,
$\gamma$ is composed of edges of length $1/n$ and $\zeta\geq2d$ so
$(\gamma_t)_{t_1\leq t <t_2}$ and, thus, $\gamma$ contains at least one
edge which is included in $\U(k)$. Noticing that for all $n\geq N$,
\[
\gamma \subset \Omega_n \subset \Omega',
\]
we obtain that this edge belongs to $\U(k) \cap\Omega' = \M'(k)$.
\end{pf*}


\section{\texorpdfstring{Covering of $\p P \cap\Omega'$ by cylinders}
{Covering of partial P cap Omega' by cylinders}}
\label{chapitre7cylinders}

From now on we only consider $n\geq N$. According to Lemma
\ref{chapitre7lemcutset}, we know that each set $\M'(k)$ for $k\in\{
0,\ldots,\lfloor hn/\zeta \rfloor-1 \}$
contains a $(\G^1_n, \G^2_n)$-cut in $\Omega_n$, thus, if we denote by~$M'(k)$
the set of the edges included in $\M'(k)$, we obtain
\[
\phi_n \leq \min\bigl\{V(M'(k)), k\in\{0,\ldots,\lfloor hn/\zeta
\rfloor-1 \} \bigr\}.
\]
However, we do not have estimates on $V(M'(k))$ that allow us to
control~$\phi_n$ using only the previous inequality. What we can use are
the upper large deviations for the maximal flow from the top to
the bottom of a cylinder (Theorem~\ref{thmdevsupphi}). In
this section, we will transform our family of cuts $(M'(k))$ by replacing
a huge part of the edges in each $\M'(k)$ by the edges of minimal
cutsets in
cylinders.

We denote by $H_i, i=1,\ldots,\N$, the intersection of the faces of $\p P$
with~$\Omega'$. For each $i =1,\ldots,\N$, we
denote by $v_i$ the exterior normal unit vector to~$P$ along $H_i$. We will
cover $\p P \cap\Omega'$ by cylinders except a surface of $\mathcal
{H}^{d-1}$ measure
controlled by a parameter $\eps$. To explain the construction of a
cutset we
will do with a huge number of cylinders, we present first the
simpler construction of a cutset using one cylinder. Let $R$
be a hyperrectangle that is included in $H_j$ for a $j\in\{1,\ldots,\N\}$
and let $B$ be the cylinder defined by
\[
B = \{x+t v_j \mid  x\in R, t\in[0,h] \},
\]
where $h\leq h_0$ is the same parameter as previously. The cylinder $B$ is
built on $\p P \cap\Omega'$, in $\RR^d \setminus\Po$. We recall that $h_0
= d(\p P, \G^1 \cup\G^2) >0$, so we know that $d(B, \G^1\cup\G^2)
>0$. We
denote by $E_a$ the set of the edges included in
\[
\E_a = \{ x+tv_j \mid  x\in R, d(x,\p R) < \zeta/n, t\in
[0,h]\}.
\]
The set $\E_a$ is a neighborhood in $B$ of the ``vertical'' faces of $B$,
that is, the faces of $B$ that are
collinear to $v_j$. We denote by $E_b$ a set of edges in $B$ that cuts the
top $R + h v_j$ from the bottom $R$ of $B$. Let $M'(k)$ be the set of the
edges included in $\M'(k)$, for a $k\in\{0,\ldots,\lfloor hn/\zeta
\rfloor
-1 \} $. Let $B'$ be the thinner cylinder
\[
B' = \{x+t v_j \mid  x\in R, d(x, \p R) \geq\zeta/n,
t\in[0,h] \}.
\]
Thus, for all $k\in\{0,\ldots,\lfloor hn/\zeta \rfloor
-1 \}$, the set of edges
\[
\bigl( M'(k) \cap(\RR^d \setminus B' )\bigr) \cup E_a \cup E_b
\]
cuts $\G^1_n$ from $\G^2_n$ in $\Omega_n$. Indeed, the set of edges
$M'(k)$ is already a cut between $\G^1_n$ and $\G^2_n$ in $\Omega_n$.
We remove from it the edges that are inside~%
$B'$ which is in the interior of $B$ and we add to it a cutset $E_b$ from
the top to the bottom of $B$ and the set of edges $E_a$ that glue together
$E_b$ and $M'(k) \cap(\RR^d \setminus B' ) $. This property is
illustrated in the Figure \ref{chapitre7constrcyl}.
\begin{rem}
In this figure, we have
represented $E_b$ as a surface (so a path in dimension $2$) that separates
the top from the bottom of the cylinder to illustrate
the fact that $E_b$ cuts all discrete paths from the bottom to the top of
$B$. Actually, we can mention that it is possible to define an object which
could be the dual of an edge in dimension $d\geq2$ (as a
generalization of
%
\begin{figure}

\includegraphics{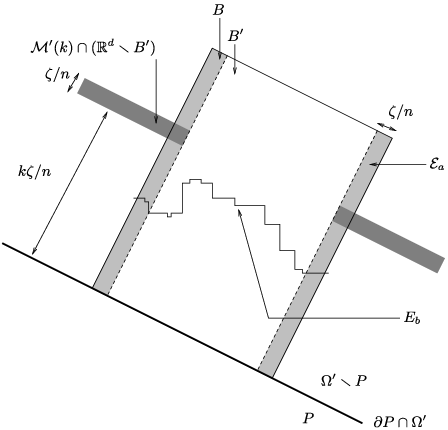}

\caption{Construction of a $(\G^1_n, \G^2_n)$-cut in $\Omega_n$ using a cutset
in a cylinder.}
\label{chapitre7constrcyl}
\end{figure}
the dual of a planar graph). This object is a plaquette,
that is, a hypersquare of sidelength $1/n$ that is orthogonal to the
edge and
cuts it in its middle and whose sides are parallel to the hyperplanes of
the axis. Then the dual of a~cutset is a hypersurface of plaquettes, thus,
Figure~\ref{chapitre7constrcyl} is somehow intuitive.
\end{rem}

\begin{figure}[t]

\includegraphics{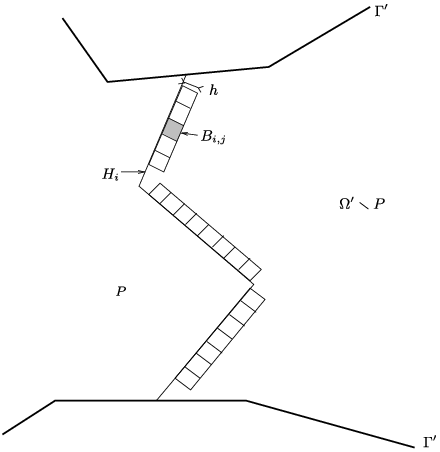}

\caption{Covering of $\p P \cap\Omega'$ by cylinders.}
\label{chapitre7polyhedre}
\end{figure}

We use exactly the same construction, but with a large number of
cylinders, that
will almost cover $\p P \cap\Omega'$. We consider a fixed $\eps
>0$. There exists a~$l$ sufficiently\vadjust{\goodbreak} small
(depending on $F$, $P$ and $\eps$) such that there
exists a finite collection $(R_{i,j}, i=1,\ldots,\N, j=1,\ldots,N_i)$ of
hypersquares of side $l$ of disjoint interiors satisfying $R_{i,j}
\subset
H_i$ for all $i\in\{1,\ldots,\N\}$ and $j\in\{1,\ldots,N_i\}$, and for all
$i\in\{1,\ldots,\N\}$,
\begin{eqnarray*}
&&\{ x\in H_i \mid  d(x, \p H_i) \geq\eps\mathcal{H}^{d-2}(\p H_i)^{-1}
\N^{-1}
\} \\
&&\qquad\subset \bigcup_{j=1}^{N_i} R_{i,j}
\subset \{ x\in H_i \mid
d(x, \p H_i) \geq\eps\mathcal{H}^{d-2}(\p H_i)^{-1} \N^{-1} 2^{-1}\}.
\end{eqnarray*}
We immediately obtain that
\[
\mathcal{H}^{d-1} \Biggl( (\p P \cap\Omega') \Bigm\backslash\bigcup_{i=1}^{\N}
\bigcup_{j=1}^{N_i} R_{i,j} \Biggr) \leq \eps.
\]
We remark that
\[
\int_{\p P \cap\Omega'} \nu(v_P(x)) \,d\mathcal{H}^{d-1}(x) \geq \sum
_{i=1}^{\N}
N_i l^{d-1}\nu(v_i),
\]
so that
\[
\PP[\phi_n\geq\lambda n^{d-1}] \leq \PP\Biggl[\phi_n \geq
(1+s/2) n^{d-1} \sum_{i=1}^{\N} N_i l^{d-1}\nu(v_i)\Biggr].
\]
Let $h<h_0$. For all $i\in\{1,\ldots,\N\}$ and $j\in\{1,\ldots,N_{i}\}$, we define
\[
B_{i,j} = \{ x + tv_i \mid  x\in R_{i,j}, t\in[0,h] \}.
\]
Since all the $B_{i,j}$ are at strictly positive distance of $\p H_i$,
there exists a~positive $h_1$ such that for all $h<h_1$, the cylinders
$B_{i,j}$ have pairwise disjoint interiors. We thus consider $h< \min
(h_0, h_1)$
(see, e.g., Figure \ref{chapitre7polyhedre}).
%
At this point, we could define a neighborhood of the vertical faces of
each cylinder $B_{i,j}$, and do the same construction as in the previous
example with one cylinder. Actually, we need to choose a little bit more
carefully the sets of edges we define along the vertical faces of the
cylinders. We will not consider only each cylinder $B_{i,j}$, but also
thinner versions of these cylinders of the type
\[
B_{i,j}(k) = \{ x+t v_j \mid  x\in R_{i,j}, d(x, \p R_{i,j}) >
k \zeta/n, t\in[0,h]\}
\]
for different values of $k$. We will then
consider the edges included in a neighborhood of the vertical faces of
each $B_{i,j }(k)$ [see the set $\W_{i,j}(k)$ above] and choose $k$ to
minimize the capacity of the union over
$i$ and $j$ of these edges. The reason why we need this optimization is
also the reason why we built a~family $(M'(k))$ of cutsets and not only one
cutset from $\G^1_n$ to $\G^2_n$ in $\Omega_n$; we will try to explain
it in
Remark \ref{chapitre7explication}.

\begin{figure}[t]

\includegraphics{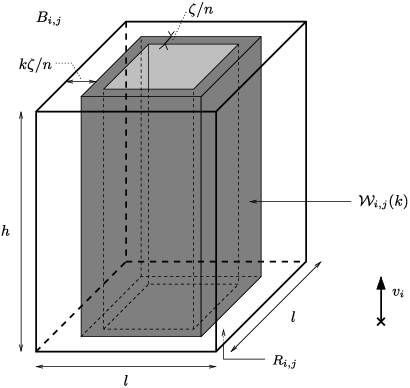}

\caption{The set $\mathcal{W}_{i,j}(k)$.}
\label{chapitre7ensW}
\end{figure}

\begin{figure}

\includegraphics{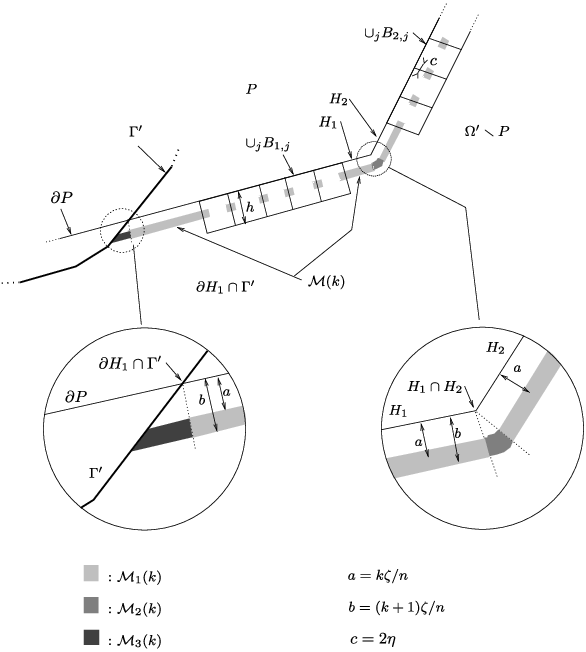}

\caption{The set $\mathcal{M}(k)$.}
\label{chapitre7ensM}
\end{figure}

Here are the precise definitions of the sets of edges.
We still consider the same constants $\zeta$ bigger than $2d$ and
$h<\min(h_0,h_1)$.
We define another positive constant $\eta$ that we will choose later
(depending on $P$, $s$ and $\Omega$). For~$i$ in $\{1,\ldots,\N\}$
and $j$ in $\{1,\ldots,N_i\}$ we recall the definition of $B_{i,j}$:
\[
B_{i,j} = \{ x+tv_i \mid  x\in R_{i,j}, t\in[0,h] \}
\]
and we define the following subsets of $\RR^d$:
%
\begin{eqnarray*}
&&\hspace*{39.4pt}B'_{i,j} = \{ x+tv_i \mid  x\in R_{i,j}, d(x, \p R_{i,j})
> \eta, t\in[0,h] \},
\\
&&\forall k \in\{0,\ldots,\lfloor\eta n / \zeta-1 \rfloor\}\\
&&\qquad\mathcal{W} _{i,j} (k) = \{ x\in B_{i,j} \mid  k\zeta/n \leq
d_2(x, \p R_{i,j} + \RR v_i) < (k+1)\zeta/n \} ,
\\
&&\forall k \in \{0,\ldots, \lfloor h n \kappa/ \zeta-1 \rfloor
\}\\
&&\qquad\hspace*{8.65pt}\M(k) = \M'(k) \Bigm\backslash\biggl( \bigcup_{i,j} B'_{i,j}\biggr)
\end{eqnarray*}
(see Figures \ref{chapitre7ensW} 
and \ref{chapitre7ensM}).
%
We denote by $W_{i,j}(k)$ the set of the edges included in
$\mathcal{W}_{i,j} (k)$ and we define $W(k)= \bigcup_{i,j} W_{i,j} (k)$. We
also denote by $M(k)$ the edges included in $\mathcal{M}(k)$. Exactly
as in
%
%
the construction of a cutset with one cylinder, we obtain a cutset that is
built with cutsets in each cylinders $B_{i,j}$. Indeed, if we denote by
$E_{i,j}$ a set of edges that is a
cutset from the top to the bottom of $B_{i,j}$ (oriented toward the
direction given by $v_i$), then for each $k_1 \in\{ 0,\ldots,\lfloor\eta
n /\zeta-1
\rfloor\}$ and $k_2 \in\{0,\ldots,\lfloor hn /\zeta-1 \rfloor\}$,
the set of edges
\[
\mathop{\bigcup_{i=1,\ldots,\N}}_{j=1,\ldots,N_i}
E_{i,j} \cup W(k_1) \cup M(k_2)
\]
contains a cutset from $\G^1_n$ to $\G^2_n$ in $\Omega_n$.
We deduce that
%
\begin{equation}
\label{chapitre7eqcontrole}
\phi_n \leq \sum_{i,j} \phi_{B_{i,j}} + \min_{k_1} V(W(k_1)) +
\min_{k_2} V( M(k_2)).
\end{equation}


\section{Control of the cardinality of the sets of edges $W$ and $M$}

For the sake of clarity, we do
not recall the sets in which the parameters
take its values; we always assume that they are the following: $i \in
\{1,\ldots,\N\}$, $j\in\{1,\ldots,N_i \}$, $k_1 \in\{0,\ldots,\lfloor\eta n /
\zeta-1 \rfloor\}$ and $k_2 \in \{ 0,\ldots,
\lfloor h n / \zeta-1 \rfloor\} $.
We have to evaluate the number of edges in the sets $W(k_1)$ and
$M(k_2)$ to control the terms $\min_{k_1} V(W(k_1))$ and $\min_{k_2} V(
M(k_2))$ in (\ref{chapitre7eqcontrole}). There
exist constants $c_1(d, \Omega) $, $c_2(P,d,\Omega)$ such that
\[
\card W(k_1) \leq c_1 \frac{\mathcal{H}^{d-1}(\p P \cap\Omega
')}{l^{d-1}} \zeta l^{d-2}
h n^{d-1} \leq c_2 l^{-1} h n^{d-1}.
\]
The cardinality of $M(k_2)$ is a little bit more complicated to
control. We
will divide $M(k)$ [resp., $\mathcal{M}(k)$] into three parts: $M(k)
\subset M_1(k) \cup M_2(k) \cup M_3(k)$ [resp., $\mathcal{M}(k)
\subset\mathcal{M}_1(k) \cup\mathcal{M}_2(k) \cup\mathcal{M}_3(k) $],
that are represented in Figure \ref{chapitre7ensM}.

We define $R'_{i,j} = \{x\in R_{i,j} \mid  d(x, \p R_{i,j})>\eta\} $ which
is the basis of $B'_{i,j}$. The set $\mathcal{M}_1(k)$ is a translation of
the sets $H_i \setminus(\bigcup_{j=1}^{N_i} R'_{i,j})$ along the
direction given by $v_i$ enlarged with a thickness $\zeta/(n\kappa)$,
\[
\mathcal{M}_1(k) \subset \bigcup_{i=1}^{\N} \Biggl\{ x + tv_i \mid  x\in H_i
\Bigm\backslash\Biggl(\bigcup_{j=1}^{N_i} R'_{i,j}\Biggr), t \in[k\zeta/ n,
(k+1)\zeta/n [ \Biggr\}.
\]
Here we have an inclusion and not an equality because $\mathcal
{M}_1(k)$ can
be a truncated version of this set (truncated at the junction
between the translates of two different faces). Since we know that
\[
\mathcal{H}^{d-1} \Biggl( (\p P \cap\Omega') \Bigm\backslash\bigcup_{i=1}^{\N}
\bigcup_{j=1}^{N_i} R_{i,j} \Biggr) \leq \eps
\]
and
\[
\mathcal{H}^{d-1} \Biggl( \bigcup_{i=1}^{\N} \bigcup_{j=1}^{N_i}( R_{i,j}
\setminus R'_{i,j}) \Biggr) \leq \frac{\mathcal{H}^{d-1}(\p P \cap\Omega
')}{l^{d-1}} l^{d-2} \eta =
\mathcal{H}^{d-1}(\p P \cap\Omega') l^{-1} \eta,
\]
we have the following bound on the cardinality of $M_1(k)$:
\[
\card(M_1(k)) \leq c_3 (\eps+ l^{-1} \eta) n^{d-1}
\]
for a constant $c_3 (d,P,\Omega, \Omega')$.

The part $M_2(k)$ corresponds to the edges
included in the ``bends'' of the neighborhood of $\p P$ located around
the boundary of the faces of $\p P$ in $\Omega'$, denoted by $\mathcal{M}_2(k)$,
that is,
\[
\M_2(k) \subset \bigcup_{i,j} \bigl( \V_2 \bigl(H_i \cap H_j, (k+1) \zeta
/n \bigr) \setminus\V_2(H_i \cap H_j, k\zeta/n ) \bigr)
\]
and there exists a constant $c_4(d, P, \Omega')$ such that
\[
\card M_2(k) \leq c_4 |k \zeta/n |^{d-2} n^{d-1} \leq c_4
h^{d-2} n^{d-1}.
\]

\begin{figure}[t]

\includegraphics{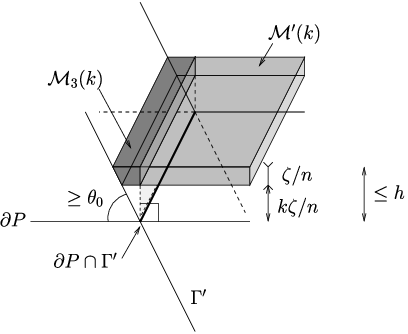}

\caption{The set $\mathcal{M}_3(k)$.}
\label{chapitre7ensM3}
\end{figure}

The last part $\M_3(k)$ corresponds to the part of $\M(k)$ that is near
the boundary $\G'$ of $\Omega'$. Indeed, $\G'$ is not orthogonal to $\p
P$,
thus, for some $k$, the set $\M(k)$ may contain edges that are not included
in
\[
\bigcup_{i=1}^{\N} \Biggl\{ x + tv_i \mid  x\in H_i
\Bigm\backslash\Biggl(\bigcup_{j=1}^{N_i} R'_{i,j}\Biggr), t \in[k\zeta/ n,
(k+1)\zeta/n [ \Biggr\},
\]
nor in
\[
\bigcup_{i,j} \bigl( \V_2 \bigl(H_i \cap H_j, (k+1) \zeta
/n \bigr) \setminus\V_2(H_i \cap H_j, k\zeta/n ) \bigr)
\]
(see Figure \ref{chapitre7ensM}). However, $\M(k) \subset\U(k)$, the
problem is to evaluate the difference of cardinality between the different
$M(k)$ due to the intersection of $\U(k)$ with $\Omega'$. We have constructed
$\Omega'$ such that $\G'$ is transverse to $\p P$ precisely to obtain\vadjust{\goodbreak} this
control. The sets $\G'$ and $\p P$ are polyhedral surfaces which are
transverse. We denote by $(\mathcal{H}_i, i\in I)$ [resp., $(\mathcal
{H}'_j,j\in J) $] the
hyperplanes that contain $\p P$ (resp., $\G'$) and by $v_i$ (resp., $v'_j$)
the exterior normal unit vector to $P$ along $\mathcal{H}_i$ (resp.,
$\Omega'$ along
$\mathcal{H}'_j$). The set $\G' \cap\p P$ is included
in the union of a finite number of intersections $\mathcal{H}_i \cap
\mathcal{H}'_j$ of
transverse hyperplanes. To each such intersection $\mathcal{H}_i \cap
\mathcal{H}'_j$, we can
associate the angles between $v_i$ and $v'_j$ and between $v_i$ and
$-v'_j$ in the plane of dimension $2$ spanned by $v_i$ and $v_j'$. Each
such angle is strictly positive because $\mathcal{H}_i$ is transverse
to $\mathcal{H}'_j$ and
so the minimum $\theta_0$ over the finite number of defined angles is
strictly positive. This $\theta_0$ and the measure $\mathcal
{H}^{d-2}(\p P \cap
\G')$ give to us a control on the volume of
$\M_3(k)$ and, thus, on $\card(M_3(k))$, as soon as these sets belong
to a
neighborhood of $\p P \cap\G'$ (see Figure \ref{chapitre7ensM3}).
%
Thus, there exist $h_2 (\Omega', P) >0$ and a
constant $c_5(d, P,\Omega, \Omega')$ such that for all $h\leq h_2$,
\[
\card(M_{3}) (k) = c_5 h n^{d-1}.
\]
We conclude that there exists a positive constant $c_6(d, P, \Omega,
\Omega')$ such that
\[
\card M(k) \leq c_6(\eps+ l^{-1} \eta+ h^{d-2} + h ) n^{d-1}.
\]


\section{Calibration of the constants}
\label{chapitre7calibrateupper}

We remark that the sets $W(k)$ [resp., the sets $M(k)$] are pairwise
disjoint for different $k$. Then we obtain that
\begin{eqnarray*}
&&\PP[\phi_n \geq \lambda n^{d-1}] \\
&&\qquad\leq \PP\Biggl[\phi_n \geq(1+s/2)
n^{d-1} \sum_{i=1}^{\N} N_i l^{d-1} \nu(v_i) \Biggr]\\
&&\qquad\leq \PP\Biggl[\sum_{i=1}^{\N} \sum_{j=1}^{N_i} \phi_{B_{i,j}} \geq
(1+s/4) n^{d-1} \sum_{i=1}^{\N} N_i l^{d-1} \nu(v_i) \Biggr]\\
&&\qquad\quad{} + \PP\Biggl[\min_{k_1} V(W(k_1)) \geq(s/8) n^{d-1}
\sum_{i=1}^{\N} N_i l^{d-1} \nu(v_i) \Biggr]\\
&&\qquad\quad{} + \PP\Biggl[\min_{k_2} V(M(k_2)) \geq(s/8) n^{d-1}
\sum_{i=1}^{\N} N_i l^{d-1} \nu(v_i) \Biggr]\\
&&\qquad\leq \sum_{i=1}^{\N} \sum_{j=1}^{N_i} \Bigl( \max_{i,j} \PP
[\phi_{B_{i,j}} \geq l^{d-1} \nu(v_i) (1+s/4) n^{d-1}] \Bigr)\\
&&\qquad\quad{} + \PP\Biggl[\sum_{i=1}^{c_2 l^{-1} h n^{d-1}} t(e_i) \geq(s/8) n^{d-1}
\sum_{i=1}^{\N} N_i l^{d-1} \nu(v_i)
\Biggr]^{\lfloor\eta n /\zeta\rfloor}\\
&&\qquad\quad{} + \PP\Biggl[\sum_{i=1}^{ c_6(\eps + l^{-1} \eta+ h^{d-2} + h) n^{d-1}
} t(e_i) \geq(s/8) n^{d-1} \sum_{i=1}^{\N} N_i l^{d-1} \nu(v_i)
\Biggr]^{2 \lfloor h n /\zeta\rfloor}.
\end{eqnarray*}
The terms
\[
\PP[\phi_{B_{i,j}} \geq l^{d-1} \nu(v_i)(1+s/4) n^{d-1}]
\]
have already been studied in \cite{Theretuppertau} (we recalled it as
Theorem \ref{thmdevsupphi} in this paper).

It remains to study two terms of the type
\[
\mathcal{P}(n) = \PP\Biggl( \sum_{i=1}^{\alpha n^{d-1}} t(e_i) \geq\beta
n^{d-1} \Biggr).
\]
As soon as $\beta> \alpha\EE(t) $ and the law of the capacity of the
edges admits an exponential moment, the Cram\'{e}r theorem in $\RR$
allows us
to affirm that
\[
\limsup_{n\rightarrow\infty} \frac{1}{n^{d-1}} \log\mathcal{P}(n)
< 0.
\]
Moreover, for all
\[
\eps \leq \eps_0 = \frac{1}{2\nu_{\max}} \int_{\mathcal{P}\cap\Omega
'} \nu(v_P(x))
\,d\mathcal{H}^{d-1}(x) ,
\]
we have
\begin{eqnarray*}
\sum_{i=1}^{\N} N_i l^{d-1} \nu(v_i)& \geq& \int_{\p P \cap\Omega'}
\nu(v_{P}(x)) \,d\mathcal{H}^{d-1} (x) - \eps\nu_{\max}\\
& \geq& \frac{1}{2} \int_{\p P \cap\Omega'} \nu(v_{P}(x)) \,d\mathcal
{H}^{d-1} (x)\\
& \geq& \frac{\nu_{\min}}{2} \mathcal{H}^{d-1}(\p P \cap\Omega').
\end{eqnarray*}
Thus, for all $\eps< \eps_0$ and $h<\min(h_0, h_1, h_2)$, if the constants
satisfy the two following conditions:
%
\begin{equation}
\label{chapitre7cond4}
c_2 l^{-1} h < \mathcal{H}^{d-1}(\p P \cap\Omega') \nu_{\min} \EE
(t(e)) s/16
\end{equation}
and
%
\begin{equation}
\label{chapitre7cond5}
c_6(\eps+ l^{-1} \eta+ h^{d-2} + h) < \mathcal{H}^{d-1}(\p P
\cap\Omega') \nu_{\min} \EE(t(e)) s/16,
\end{equation}
thanks to Theorem \ref{thmdevsupphi} and the Cram\'{e}r theorem in $\RR
$, we obtain that
\[
\limsup_{n\rightarrow\infty} \frac{1}{n^{d}}\log\PP[\phi_n \geq
\lambda n^{d-1}] < 0
\]
and Theorem \ref{chapitre7devsup} is proved. We claim that it is
possible to
choose the constants such that conditions (\ref{chapitre7cond4}) and
(\ref{chapitre7cond5}) are satisfied. Indeed, we first choose $\eps
<\eps_0$ such that
\[
\eps < \frac{1}{4} \frac{\mathcal{H}^{d-1} (\p P \cap\Omega) \nu_{\min
} \EE(t(e)) s}{16
c_6}.
\]
To this fixed $\eps$ corresponds a $l$. Knowing $\eps$ and $l$, we
choose $h \leq\min(h_0, h_1, h_2)$ and $\eta$
such that
\[
\max(h, h^{d-2},l^{-1} h, l^{-1} \eta) < \frac{1}{4} \frac{
\mathcal{H}^{d-1} (\p P \cap\Omega') \nu_{\min}
\EE(t(e)) s}{16 \max(c_2, c_6)}.
\]
The proof of Theorem \ref{chapitre7devsup} is complete.
\begin{rem}
\label{chapitre7explication}
We try here to explain why we built several sets $W(k_1)$ and $M(k_2)$, and
not only one couple of such sets, that would have been sufficient to
construct a cutset from
$\G^1_n$ to $\G^2_n$ in $\Omega_n$.
To use estimates of upper large deviations of maximal flows in cylinder we
already know, we want to compare $\phi_n$ with $\sum_{i,j}
\phi_{B_{i,j}}$. Heuristically, to construct a $(\G^1_n, \G^2_n)$-cut in
$\Omega_n$ from the union of cutsets in each cylinder $B_{i,j}$, we
have to add edges to glue together the different cutsets at the common
boundary of the small cylinders and to extend these cutsets to $ (\p P
\cap\Omega_n) \setminus\bigcup_{i=1}^{\N} \bigcup_{j=1}^{N_i}
R_{i,j}$. Yet we
want to prove that the upper large deviations of $\phi_n$ are of
volume order. If we only consider one possible set $E$ of edges such that
\[
\phi_n \leq \sum_{i,j} \phi_{B_{i,j}} + V(E),
\]
we will obtain that
\begin{eqnarray*}
\PP[\phi_n \geq\lambda n^{d-1} ]& \leq &\sum_{i,j} \PP[\phi_{B_{i,j}}
\geq l^{d-1}\nu(v_i) (1+s/4 ) n^{d-1}]\\
&&{} + \PP\Biggl[V(E) \geq
n^{d-1} \sum_{i=1}^{\N} N_i l^{d-1}\nu(v_i) s/4\Biggr].
\end{eqnarray*}
We can choose such a set $E$ so that it contains less than $\delta n^{d-1}$
edges for a~small $\delta$ [e.g., $E$ is equal to $W(k_1) \cup M(k_2)$
for a
fixed couple $(k_1,k_2)$] but the probability
\[
\PP\Biggl[\sum_{i=1}^{\delta n^{d-1}} t(e_i) \geq C n^{d-1}\Biggr]
\]
does not decay exponentially fast with $n^d$ in general. To obtain this
speed of
decay, we have to make an optimization over the possible choices of the
set~$E$, that is, we choose~$E$ among a set of $C' n$ possible disjoint
sets of edges
$E_1,\ldots,E_{C'n}$; in this case, we obtain that
\[
\phi_n \leq \sum_{i,j} \phi_{B_{i,j}} +
\min_{k=1,\ldots,C'n}V(E_k)
\]
and so
%
\begin{eqnarray}
\label{chapitre7heuristique}
\PP[\phi_n \geq\lambda n^{d-1} ]& \leq & \sum_{i,j} \PP[\phi_{B_{i,j}}
\geq l^{d-1}\nu(v_i) (1+s/4 ) n^{d-1}] \nonumber\\[-8pt]\\[-8pt]
&&{} + \prod_{k=1}^{C'n} \PP\Biggl[V(E_k) \geq n^{d-1} \sum_{i=1}^{\N} N_i
l^{d-1}\nu(v_i) s/4 \Biggr].\nonumber
\end{eqnarray}
It is then sufficient to prove that for all $k$, $ \PP[V(E_k) \geq C''
n^{d-1}] $ decays exponentially fast with $n^{d-1}$ to conclude that
the last
term in (\ref{chapitre7heuristique}) decays exponentially fast with
$n^d$. Theorem \ref{thmdevsupphi} gives a control on the terms
\[
\PP[\phi_{B_{i,j}} \geq l^{d-1}\nu(v_i) (1+s/4 ) n^{d-1}].
\]
The conclusion is that to obtain the volume order of the upper large
deviations, the optimization over the different possible values of $k_1$
and $k_2$ is really important, even if it is not needed if we only want to
prove that $\PP(\phi_n \geq\lambda n^{d-1})$ goes to zero when $n$ goes
to infinity.
\end{rem}


%
\printaddresses


\begin{thebibliography}{17}

\bibitem{Bollobas}
\begin{bbook}[mr]
\bauthor{\bsnm{Bollob{\'a}s},~\bfnm{B{\'e}la}\binits{B.}}
(\byear{1979}).
\btitle{Graph Theory: An Introductory Course}.
\bseries{Graduate Texts in Mathematics}
\bvolume{63}.
\bpublisher{Springer}, \baddress{New York}.
\bid{mr={0536131}}
\end{bbook}
\endbibitem

\bibitem{CerfStFlour}
\begin{bbook}[mr]
\bauthor{\bsnm{Cerf},~\bfnm{R.}\binits{R.}}
(\byear{2006}).
\btitle{The {W}ulff Crystal in {I}sing and Percolation Models}.
\bseries{Lecture Notes in Math.}
\bvolume{1878}.
\bpublisher{Springer}, \baddress{Berlin}.
\bid{mr={2241754}}
\end{bbook}
\endbibitem

\bibitem{CerfTheret09geob}
\begin{bmisc}[auto:STB|2010-11-18|09:18:59]
\bauthor{\bsnm{Cerf},~\bfnm{Rapha{\"e}l}\binits{R.}} \AND
  \bauthor{\bsnm{Th{\'e}ret},~\bfnm{Marie}\binits{M.}}
(\byear{2011}).
\bhowpublished{Law of large numbers for the maximal flow through a domain of
  $\mathbb{R}^d$ in first passage percolation. \textit{Trans. Amer. Math.
  Soc.} \textbf{363} 3665--3702}.
\end{bmisc}
\endbibitem

\bibitem{CerfTheret09infb}
\begin{bmisc}[auto:STB|2010-11-18|09:18:59]
\bauthor{\bsnm{Cerf},~\bfnm{Rapha{\"e}l}\binits{R.}} \AND
  \bauthor{\bsnm{Th{\'e}ret},~\bfnm{Marie}\binits{M.}}
(\byear{2011}).
\bhowpublished{Lower large deviations for the maximal flow through a domain of
  $\mathbb{R}^d$ in first passage percolation. \textit{Probab. Theory Related
  Fields}. To appear. Available at
  \href{http://arxiv.org/abs/0907.5501}{arxiv.org/abs/0907.5501}%
}.
\end{bmisc}
\endbibitem

\bibitem{FED}
\begin{bbook}[mr]
\bauthor{\bsnm{Federer},~\bfnm{Herbert}\binits{H.}}
(\byear{1969}).
\btitle{Geometric Measure Theory}.
\bseries{Die Grundlehren der Mathematischen Wissenschaften}
\bvolume{153}.
\bpublisher{Springer}, \baddress{New York}.
\bid{mr={0257325}}
\end{bbook}
\endbibitem

\bibitem{Garet2}
\begin{barticle}[mr]
\bauthor{\bsnm{Garet},~\bfnm{Olivier}\binits{O.}}
(\byear{2009}).
\btitle{Capacitive flows on a 2{D} random set}.
\bjournal{Ann. Appl. Probab.}
\bvolume{19}
\bpages{641--660}.
\bid{doi={10.1214/08-AAP556}, mr={2521883}}
\end{barticle}
\endbibitem

\bibitem{GrimmettKesten84}
\begin{barticle}[mr]
\bauthor{\bsnm{Grimmett},~\bfnm{Geoffrey}\binits{G.}} \AND
  \bauthor{\bsnm{Kesten},~\bfnm{Harry}\binits{H.}}
(\byear{1984}).
\btitle{First-passage percolation, network flows and electrical resistances}.
\bjournal{Z. Wahrsch. Verw. Gebiete}
\bvolume{66}
\bpages{335--366}.
\bid{doi={10.1007/BF00533701}, mr={0751574}}
\end{barticle}
\endbibitem

\bibitem{KestenStFlour}
\begin{bincollection}[mr]
\bauthor{\bsnm{Kesten},~\bfnm{Harry}\binits{H.}}
(\byear{1986}).
\btitle{Aspects of first passage percolation}.
In \bbooktitle{\'{E}cole D'\'et\'e de Probabilit\'es de {S}aint-{F}lour,
  {XIV}---1984}.
\bseries{Lecture Notes in Math.}
\bvolume{1180}
\bpages{125--264}.
\bpublisher{Springer}, \baddress{Berlin}.
\bid{mr={0876084}}
\bptnote{check year}
\end{bincollection}
\endbibitem

\bibitem{Kestenflows}
\begin{barticle}[mr]
\bauthor{\bsnm{Kesten},~\bfnm{Harry}\binits{H.}}
(\byear{1987}).
\btitle{Surfaces with minimal random weights and maximal flows: A~%
  higher-dimensional version of first-passage percolation}.
\bjournal{Illinois J. Math.}
\bvolume{31}
\bpages{99--166}.
\bid{mr={0869483}}
\end{barticle}
\endbibitem

\bibitem{LA}
\begin{bbook}[mr]
\bauthor{\bsnm{Lang},~\bfnm{Serge}\binits{S.}}
(\byear{1985}).
\btitle{Differential Manifolds}, \bedition{2nd} ed.
\bpublisher{Springer}, \baddress{New York}.
\bid{mr={0772023}}
\end{bbook}
\endbibitem

\bibitem{RossignolTheret09}
\begin{barticle}[mr]
\bauthor{\bsnm{Rossignol},~\bfnm{Rapha{\"e}l}\binits{R.}} \AND
  \bauthor{\bsnm{Th{\'e}ret},~\bfnm{Marie}\binits{M.}}
(\byear{2010}).
\btitle{Law of large numbers for the maximal flow through tilted cylinders in
  two-dimensional first passage percolation}.
\bjournal{Stochastic Process. Appl.}
\bvolume{120}
\bpages{873--900}.
\bid{doi={10.1016/j.spa.2010.02.005}, mr={2610330}}
\end{barticle}
\endbibitem

\bibitem{RossignolTheret08b}
\begin{barticle}[auto:STB|2010-11-18|09:18:59]
\bauthor{\bsnm{Rossignol},~\bfnm{Rapha{\"e}l}\binits{R.}} \AND
  \bauthor{\bsnm{Th{\'e}ret},~\bfnm{Marie}\binits{M.}}
(\byear{2010}).
\btitle{Lower large deviations and laws of large numbers for maximal
  flows through a box in first passage percolation}.
\bjournal{Ann. Inst. H. Poincar\'{e} Probab. Statist.}
\bvolume{46}
\bpages{1093--1131}.
\bid{mr={2744888}}
\end{barticle}
\endbibitem

\bibitem{TheretUpper}
\begin{barticle}[mr]
\bauthor{\bsnm{Th{\'e}ret},~\bfnm{Marie}\binits{M.}}
(\byear{2007}).
\btitle{Upper large deviations for the maximal flow in first-passage
  percolation}.
\bjournal{Stochastic Process. Appl.}
\bvolume{117}
\bpages{1208--1233}.
\bid{doi={10.1016/j.spa.2006.12.007}, mr={2343936}}
\end{barticle}
\endbibitem

\bibitem{Theretuppertau}
\begin{bmisc}[auto:STB|2010-11-18|09:18:59]
\bauthor{\bsnm{Th{\'e}ret},~\bfnm{Marie}\binits{M.}}
(\byear{2009}).
\bhowpublished{Upper large deviations for maximal flows through a tilted
  cylinder.
\textit{ESAIM Probab. Stat.} To appear.
  Available at
  \href{http://arxiv.org/abs/0907.0614}{arxiv.org/abs/0907.0614}%
}.
\end{bmisc}
\endbibitem

\bibitem{Zhang}
\begin{barticle}[mr]
\bauthor{\bsnm{Zhang},~\bfnm{Yu}\binits{Y.}}
(\byear{2000}).
\btitle{Critical behavior for maximal flows on the cubic lattice}.
\bjournal{J. Stat. Phys.}
\bvolume{98}
\bpages{799--811}.
\bid{doi={10.1023/A:1018631726709}, mr={1749233}}
\end{barticle}
\endbibitem

\bibitem{Zhang07}
\begin{bmisc}[auto:STB|2010-11-18|09:18:59]
\bauthor{\bsnm{Zhang},~\bfnm{Yu}\binits{Y.}}
(\byear{2007}).
\bhowpublished{Limit theorems for maximum flows on a lattice. Available at
  \href{http://arxiv.org/abs/0710.4589}{arxiv.org/abs/0710.4589}%
}.
\end{bmisc}
\endbibitem

\end{thebibliography}
\end{document}